\documentclass[aap,seceqn,MSNbibl,citesort,dvips]{arximspdf}
\usepackage{graphicx}

\doi{10.1214/10-AAP731}
\volume{21}
\issue{4}
\pubyear{2011}
\firstpage{1466}
\lastpage{1492}

\makeatletter
\newcommand{\eqref}[1]{(\ref{#1})}

\newtheorem{theorem}{Theorem}[section]
\newproclaim{remark}[theorem]{Remark}
\newtheorem{lemma}[theorem]{Lemma}
\newtheorem{proposition}[theorem]{Proposition}

\newcommand{\del}{\partial}
\newcommand{\delt}{\del_t}
\newcommand{\lap}{\triangle}
\newcommand{\inv}{^{-1}}
\newcommand{\transpose}{^*}
\newcommand{\grad}{\nabla}
\newcommand{\gradt}{\grad\transpose}
\newcommand{\divergence}{\grad\cdot}
\newcommand{\curl}{\grad\times}

\renewcommand{\epsilon}{\varepsilon}

\newcommand{\lhp}{\mathbf{P}}

\newcommand{\R}{\mathbb{R}}
\newcommand{\N}{\mathbb{N}}

\newcommand{\at}{|}
\newcommand{\defeq}{\stackrel{\mathrm{def}}{=}}

\newcommand{\Chi}[1]{\chi_{\{#1\}}}
\newcommand{\varmax}{\vee}
\newcommand{\varmin}{\wedge}

\newcommand{\sidenote}[1]{}

\def\bsuffix #1{#1}

\def\bptnote#1{}

\makeatother

\begin{document}
\begin{frontmatter}

\title{A stochastic-Lagrangian approach to the Navier--Stokes
equations in domains with boundary}
\runtitle{Stochastic-Lagrangian Navier--Stokes with boundary}

\begin{aug}
\author[A]{\fnms{Peter} \snm{Constantin}\ead[label=e1]{const@math.uchicago.edu}\thanksref{aut1}}
and
\author[B]{\fnms{Gautam} \snm{Iyer}\corref{}\ead[label=e2]{gautam@math.cmu.edu}\thanksref{aut2}}
\runauthor{P. Constantin and G. Iyer}
\thankstext{aut1}{Supported in part by NSF Grant DMS-08-04380.}
\thankstext{aut2}{Supported in part by NSF Grant DMS-07-07920 and the
Center for Nonlinear Analysis.}
\affiliation{University of Chicago and Carnegie Mellon University}
\address[A]{Department of Mathematics\\
University of Chicago\\
Chicago, Illinois 60637\\
USA\\
\printead{e1}} 
\address[B]{Department of Mathematical Sciences\\
Carnegie Mellon University\\
Pittsburgh, Pennsylvania 15213\\
USA\\
\printead{e2}}
\end{aug}

\received{\smonth{3} \syear{2010}}
\revised{\smonth{8} \syear{2010}}

%
\begin{abstract}
In this paper we derive a probabilistic representation of the
deterministic $3$-dimensional Navier--Stokes equations in the presence
of spatial boundaries. The formulation in the absence of spatial
boundaries was done by the authors in [\textit{Comm.
Pure Appl. Math.} \textbf{61} (2008) 330--345]. While the formulation in
the presence of boundaries is similar in spirit, the proof is somewhat
different. One aspect highlighted by the formulation in the presence of
boundaries is the nonlocal, implicit influence of the boundary
vorticity on the interior fluid velocity.
\end{abstract}

%

\end{frontmatter}

\section{Introduction}\label{sxnIntroduction}
The (unforced) incompressible Navier--Stokes equations
%
\begin{eqnarray}
\label{eqnNavierStokes}
\del_t u + (u\cdot\grad) u - \nu\lap u +
\grad p &=& 0,\\
\label{eqnIncompressibility} \divergence u &=& 0
\end{eqnarray}
describe the evolution of the velocity field $u$ of an incompressible
fluid with kinematic viscosity $\nu> 0$ in the absence of external
forcing. Here $u = u(x,t)$ with $t \geq0$, $x \in\R^d$, $d \geq2$.
Equation \eqref{eqnIncompressibility} is the incompressibility
constraint. Unlike compressible fluids, the pressure $p$ in \eqref
{eqnNavierStokes} does not have a physical meaning\sidenote{Too bold?
Might offend someone.} and is only a Lagrange multiplier that ensures
incompressibility is preserved. While equations \eqref{eqnNavierStokes} and \eqref{eqnIncompressibility} can be formulated in
any dimension $d \geq2$, they are usually only studied in the
physically relevant dimensions $2$ or $3$. The presentation of the
Navier--Stokes equations above is in the absence of spatial boundaries;
an issue that will be discussed in detail later.

When $\nu= 0$, \eqref{eqnNavierStokes} and \eqref{eqnIncompressibility} are known as the Euler equations. These describe
the evolution of the velocity field of an (ideal) inviscid and
incompressible fluid. Formally the difference between the Euler and
Navier--Stokes equations is only the dissipative Laplacian term. Since
the Laplacian is exactly the generator a Brownian motion, one would
expect to have an exact stochastic representation of \eqref{eqnNavierStokes} and \eqref{eqnIncompressibility} which is physically
meaningful, that is, can be thought of as an appropriate average of the
inviscid dynamics and Brownian motion.

The difficulty, however, in obtaining such a representation is because
of both the nonlinearity and the nonlocality of equations \eqref{eqnNavierStokes} and \eqref{eqnIncompressibility}. In 2D, an exact
stochastic representation of \eqref{eqnNavierStokes} and \eqref{eqnIncompressibility} dates back to Chorin~\cite{bblChorin} in 1973
and was obtained using vorticity transport and the Kolmogorov
equations. In three dimensions, however, this method fails to provide
an exact representation because of the vortex stretching term.

In 3D, a variety of techniques has been used to provide exact
stochastic representations of \eqref{eqnNavierStokes} and \eqref{eqnIncompressibility}. One such technique (Le Jan and Sznitman~\cite
{bblLeJanSznitman}) uses a backward branching process in Fourier space.
This approach has been extensively studied and generalized \cite
{bblBhattacharyaORST1,bblBhattacharyaORST2,bblThomannOssiander,bblOssiander,bblWaymire}
by many authors (see also \cite{bblWaymire2}). A different and more
recent technique due to Busnello, Flandoli and Romito \cite
{bblBusnelloFlandolliRomito} (see also \cite{bblBusnello}) uses noisy
flow paths and a Girsanov transformation. A related approach in \cite
{bblConstIyerSLNS} is the stochastic-Lagrangian formulation, exact
stochastic representation of solutions to \eqref{eqnNavierStokes} and \eqref{eqnIncompressibility} which is essentially
the averaging of noisy particle trajectories and the inviscid dynamics.
Stochastic variational approaches (generalizing Arnold's \cite
{bblArnold} deterministic variational formulation for the Euler
equations) have been used by \cite{bblCiprianoCruzerio,bblEyink} and
a related approach using stochastic differential geometry can be found
in \cite{bblGliklikh}.

One common setback in all the above methods is the inability to deal
with boundary conditions. The main contribution of this paper adapts
the stochastic-Lagrangian formulation in \cite{bblConstIyerSLNS}
(where the authors only considered periodic boundary conditions or
decay at infinity) to the situation with boundaries. The usual
probabilistic techniques used to transition to domains with boundary
involve stopping the processes at the boundary. This introduces two
major problems with the techniques in \cite{bblConstIyerSLNS}. First,
stopping introduces spatial discontinuities making the proof used in
\cite{bblConstIyerSLNS} fail and a different approach is required.
Second and more interesting is the fact that merely stopping does
\textit{not} give the no-slip ($0$-Dirichlet) boundary condition as
one would expect. One needs to also create trajectories at the boundary
which essentially propagate the influence of the vorticity at the
boundary to the interior fluid velocity.

\subsection{Plan of the paper}
This paper is organized as follows. In Section~\ref{sxnSLIntro} a~%
brief introduction to the stochastic-Lagrangian formulation without
boundaries is given. In Section \ref{sxnSLWithBoundaries} we motivate
and state the stochastic-Lagrangian formulation in the presence of
boundaries (Theorem \ref{thmSLNSNoSlip}). In Section \ref
{sxnBackwardIto} we recall certain standard facts about backward It\^o
integrals which will be used in the proof of Theorem \ref
{thmSLNSNoSlip}. In Section \ref{sxnNoSlipProof} we prove Theorem \ref
{thmSLNSNoSlip}. Finally, in Section \ref{sxnVorticityTransport} we
discuss stochastic analogues of vorticity transport and inviscid
conservation laws.

\section{The stochastic-Lagrangian formulation without
boundaries}\label{sxnSLIntro}
In this section, we provide a brief description of the
stochastic-Lagrangian formulation in the absence of boundaries. For
motivation, let us first study a Lagrangian description of the Euler
equations [equations \eqref{eqnNavierStokes} and \eqref{eqnIncompressibility} with $\nu= 0$; we will usually use a superscript
of $0$ to denote quantities relating to the Euler equations]. Let $d =
2, 3$ denote the spatial dimension and $X^0_t$ be the flow defined by
%
\begin{equation}\label{eqnEulerXdef}
\dot X^0_t = u^0_t(X^0_t),
\end{equation}
with initial data $X^0_0(a) = a$, for all $a \in\R^d$. To
clarify our notation, $X^0$ is a~function of the initial data $a \in\R
^d$ and time $t \in[0, \infty)$. We usually omit the
spatial variable and use $X^0_t$ to denote $X^0( \cdot, t)$, the slice
of $X^0$ at time $t$. Time derivatives will always be denoted by a dot
or $\delt$ instead of a $t$ subscript.

One can immediately check (see, e.g., \cite{bblConstELE}) that
$u$ satisfies the incompressible Euler equations if and only if $\ddot
X^0$ is a gradient composed with $X$. By Newton's second law, this
admits the physical interpretation that the Euler equations are
equivalent to assuming that the force on individual particles is a gradient.

One would naturally expect that solutions to the Navier--Stokes
equations can be obtained similarly by adding noise to particle
trajectories and averaging. However, for noisy trajectories, an
assumption on $\ddot X^0$ will be problematic. In the incompressible
case, we can circumvent this difficultly using the Weber formula \cite
{bblWebber} [equation \eqref{eqnEulerWebber} below]. Indeed, a direct
computation (see, e.g., \cite{bblConstELE}) shows that for
divergence free $u$, the assumption that~$\ddot X^0$ is a~gradient is
equivalent to
%
\begin{equation}\label{eqnEulerWebber}
u^0_t = \lhp[ (\gradt A^0_t) (u^0_0 \circ A^0_t) ],
\end{equation}
where $\lhp$ denotes the Leray--Hodge projection \cite
{bblChorinMarsden,bblConstFoias,bblMajdaBertozzi} onto divergence free
vector fields, the notation $\gradt$ denotes the transpose of the
Jacobian and for any $t \geq0$, $A^0_t = (X^0_t)\inv$ is the spatial
inverse of the map $X^0_t$ [i.e.,\ $A^0_t(X^0_t(a)) = a$ for all $a \in
\R^d$ and $X^0_t(A^0_t(x)) = x$ for all $x\in\R^d$].

\setcounter{footnote}{2}

From this we see that the Euler equations are formally equivalent to
equations \eqref{eqnEulerXdef} and \eqref{eqnEulerWebber}. Since
these equations no longer involve second (time) derivatives of the flow
$X^0$, one can consider noisy particle trajectories without any
analytical difficulties. In fact, adding noise to \eqref{eqnEulerXdef}
and averaging out the noise in \eqref{eqnEulerWebber} gives the
equivalent formulation of the Navier--Stokes equations stated below.

\begin{theorem}[(Constantin, Iyer \cite{bblConstIyerSLNS})]\label{thmSLNS}
Let $d \in\{2,3\}$ be the spatial dimension, $\nu> 0$ represent the
kinematic viscosity and $u_0$ be a divergence free, periodic, H\"older
$2+\alpha$ function and $W$ be a $d$-dimensional Wiener process.
Consider the system
%
\begin{eqnarray}
\label{eqnXdef} dX_t &=& u_t(X_t) \,dt + \sqrt{2\nu} \,dW_t,\\
\label{eqnX0} X_0(a) &= &a \qquad \forall a \in\R^d,\\
\label{eqnAveragedWebber} u_t &=& E \lhp[ (\gradt A_t) (u_0 \circ
A_t) ],
\end{eqnarray}
where, as before, for any $t \geq0$, $A_t = X_t \inv$ denotes the
spatial inverse%
\footnote{It is well known (see, e.g., Kunita \cite
{bblKunita}) that the solution to \eqref{eqnXdef} and \eqref{eqnX0}
gives a stochastic flow of diffeomorphisms and, in particular,
guarantees the existence of the spatial inverse of $X$.}
of~$X_t$. Then $u$ is a classical solution of the Navier--Stokes
equations \eqref{eqnNavierStokes} and \eqref{eqnIncompressibility} with
initial data $u_0$ and periodic boundary conditions if and only if $u$
is a~fixed point of the system \eqref{eqnXdef}--\eqref{eqnAveragedWebber}.
\end{theorem}

\begin{remark*}
The flows $X, A$ above are now a function of the initial data $a \in\R
^d$, time $t \in[0, \infty)$ and the probability variable
$\varpi\in\Omega$. We always suppress the probability variable, use
$X_t$ to denote $X(\cdot, t)$ and omit the spatial variable when
unnecessary. The function $u$ is a deterministic function of space and
time and, as above, we use $u_t$ to denote the function $u(\cdot, t)$.
\end{remark*}

We now briefly explain the idea behind the proof of Theorem \ref
{thmSLNS} given in~\cite{bblConstIyerSLNS} and explain why this method
can not be used in the presence of spatial boundaries. Consider first
the solution of the SDE \eqref{eqnXdef} with initial data \eqref
{eqnX0}. Using the It\^o--Wentzel formula \cite{bblKunita}, {Theorem
4.4.5}, one can show that any (spatially regular) process $\theta$
which is constant along trajectories of~$X$ satisfies the SPDE
%
\begin{equation}\label{eqnTheta}
d\theta_t + (u_t \cdot\grad) \theta_t\, dt - \nu\lap\theta_t\, dt +
\sqrt{2\nu} \grad\theta_t \,dW_t = 0.
\end{equation}
Since the process $A$ (which, as before, is defined to be the spatial
inverse of~$X$) is constant along trajectories of $X$, the process
$\theta$ defined by
%
\begin{equation}\label{eqnThetaDef}
\theta_t = \theta_0 \circ A_t
\end{equation}
is constant along trajectories of $X$. Thus, if $\theta_0$ is regular
enough ($C^2$), then~$\theta$ satisfies SPDE \eqref{eqnTheta}. Now,
if $u$ is deterministic, taking expected values of~\eqref{eqnTheta} we
see that $\bar\theta_t = E\theta_0 \circ A_t$ satisfies
%
\begin{equation}\label{eqnBarTheta}
\del_t \bar\theta_t + (u_t \cdot\grad) \bar\theta_t - \nu\lap
\bar\theta_t = 0
\end{equation}
with initial condition $\bar\theta|_{t = 0} = \theta_0$.

\begin{remark*}
Note that when $\nu= 0$, $A$ is deterministic so $\bar\theta= E
\theta= \theta$. Further, equation \eqref{eqnTheta} reduces to the
transport equation for which writing the solution as $\theta_t =
\theta_0 \circ A_t$ is exactly the method of characteristics. When
$\nu> 0$, the above procedure is an elegant generalization, termed as
the ``method of random characteristics'' (see \cite
{bblConstIyerSLNS,bblRozovskiBook,bblIyerThesis} for further information).
\end{remark*}

Once explicit equations for $A$ and $u_0 \circ A$ have been
established, a direct computation using It\^o's formula shows that $u$
given by \eqref{eqnAveragedWebber} satisfies the Navier--Stokes
equations \eqref{eqnNavierStokes} and \eqref{eqnIncompressibility}. This
was the proof used in \cite{bblConstIyerSLNS}.

\begin{remark*}
This point of view also yields a natural understanding of generalized
relative entropies \cite
{bblMichelMischlerPerthame1,bblMichelMischlerPerthame2,bblConstEntropies,bblConstIyerSEntropy}.
Eyink's recent work \cite{bblEyinkSL} adapted this framework to
magnetohydrodynamics and related equations by using the analogous Weber
formula \cite{bblKuznetsovRuban,bblRuban}. We also mention that Zhang
\cite{bblZhang} considered a backward analogue and provided short
elegant proofs to classical existence results to \eqref{eqnNavierStokes} and \eqref{eqnIncompressibility}.
\end{remark*}

\section{The formulation for domains with boundary}\label{sxnSLWithBoundaries}

In this section we describe how \eqref{eqnXdef}--\eqref
{eqnAveragedWebber} can be reformulated in the presence of boundaries.
We begin by describing the difficulty in using the techniques from
\cite{bblConstIyerSLNS} described in Section \ref{sxnSLIntro}.

Let $D \subset\R^d$ be a domain with Lipschitz boundary. Even if we
insist $u = 0$ on the boundary of $D$, we note that the noise in \eqref
{eqnXdef} is independent of space and thus, insensitive to the presence
of the boundary. Consequently, some trajectories of the stochastic flow
$X$ will leave the domain $D$ and for any $t > 0$, the map $X_t$ will
(surely) not be spatially invertible. This renders~\eqref{eqnThetaDef}
meaningless.

\begin{figure}

\includegraphics{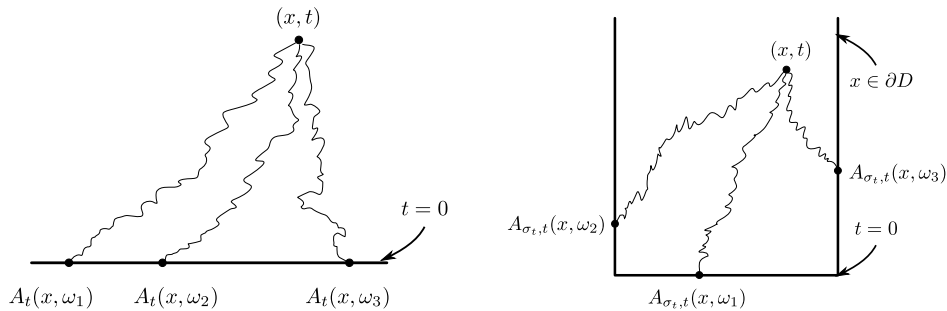}%
\vspace*{-3pt}
\caption{Three sample realizations of $A$ without boundaries \textup{(left)}
and with boundaries \textup{(right)}.}\label{fgrBoundaries}
\vspace*{-5pt}
\end{figure}

In the absence of spatial boundaries, equation \eqref{eqnThetaDef}
dictates that $\bar\theta(x,t)$ is determined by averaging the
initial data over all trajectories of $X$ which reach $x$ at time $t$.
In the presence of boundaries, one must additionally average the
boundary value of all trajectories reaching $(x,t)$, starting on $\del
D$ at any intermediate time (Figure \ref{fgrBoundaries}). As we will
see later, this means the analogue of \eqref{eqnThetaDef} in the
presence of spatial boundaries is a spatially discontinuous process.
This renders \eqref{eqnTheta} meaningless, giving a second obstruction
to using the methods of \cite{bblConstIyerSLNS} in the presence of
boundaries.

While the method of random characteristics has the above inherent
difficulties in the presence of spatial boundaries, equation \eqref
{eqnBarTheta} is exactly the Kolmogorov Backward equation (\cite
{bblOksendal}, {Section 8.1)}. In this case, an expected value
representation in the presence of boundaries is well known. More
generally, the Feynman--Kac (\cite{bblOksendal}, {Section 8.2}) formula, at
least for linear equations with a potential term, has been successfully
used in this situation. A certain version of this method (Section \ref
{sxnFeynmanKacBdd}), without making the usual time reversal
substitution, is essentially the same as the method of random
characteristics. It is this version that will yield the natural
generalization of~\eqref{eqnXdef}--\eqref{eqnAveragedWebber} in
domains with boundary (Theorem \ref{thmSLNSNoSlip}). Before turning to
the Navier--Stokes equations, we provide a brief discussion on the
relation between the Feynman--Kac formula and the method of random
characteristics.\looseness=-1

\subsection{The Feynman--Kac formula and the method of random
characteristics}\label{sxnFeynmanKacBdd}
Both the Feynman--Kac formula and the method of random characteristics
have their own advantages and disadvantages: The method of random
characteristics only involves forward SDE's and obtains the solution of
\eqref{eqnBarTheta} at time $t$ with only the knowledge of the initial
data and ``$X$ at time $t$'' (or more precisely, the solution at time
$t$ of the equation \eqref{eqnXdef} with initial data specified at
time $0$). However, this method involves computing the spatial inverse
of $X$, which analytically and numerically involves an additional step.

On the other hand, to compute the solution of \eqref{eqnBarTheta} at
time $t$ via the probabilistic representation using the Kolmogorov
backward equation (or equivalently, the Feynman--Kac formula with a $0$
potential term) when $u$ is time dependent involves \textit{backward}
SDE's and further requires the knowledge of the solution to \eqref
{eqnXdef} with initial conditions specified at all times $s \leq t$.
However, this does not require computation of spatial inverses and,
more importantly, yields the correct formulation in the presence of
spatial boundaries.

Now, to see the relation between the method of random characteristics
and the Feynman--Kac formula, we rewrite \eqref{eqnXdef} in integral
form and keep track of solutions starting at all times $s \geq0$. For
any $s \geq0$, we define the process $\{ X_{s,t} \}_{t \geq s}$ to be
the flow defined by
%
\begin{equation}\label{eqnXst}
X_{s,t}(x) = x + \int_s^t u_r \circ X_{s,r}(x) \,dr + \sqrt{2\nu} (
W_t - W_s ).
\end{equation}

Now, as always, we let $A_{s,t} \!=\! X_{s,t}\inv$. Then formally
composing \eqref{eqnXst} with~$A_{s,t}$ and using the semigroup
property $X_{s,t} \circ X_{r,s} = X_{r,t}$ gives the self-contained
backward equation for $A_{s,t}$
%
\begin{equation}\label{eqnAst}
A_{s,t}(x) = x - \int_s^t u_r \circ A_{r,t}(x)\, dr - \sqrt{2\nu} (
W_t - W_s ).
\end{equation}
Now \eqref{eqnThetaDef} can be written as
%
\begin{equation}\label{eqnThetaDefA}
\theta_t = \theta_0 \circ A_{0,t}
\end{equation}
and using the semigroup property $A_{r,s} \circ A_{s,t} = A_{r,t}$ we
see that
%
\begin{equation}\label{eqnThetaDefAst}
\theta_t = \theta_s \circ A_{s,t}.
\end{equation}
%

This formal calculation leads to a natural generalization of \eqref
{eqnThetaDef} in the presence of boundaries. As before, let $D \subset
\R^d$ be a domain with Lipschitz boundary and assume, for now, that
$u$ is a Lipschitz function defined on \textit{all} of $\R^d$. Let $A_{s,t}$ be the flow defined by \eqref{eqnAst} and for $x
\in D$, we define the \textit{backward exit time%
} $\sigma_t(x)$ by
%
\begin{equation}\label{eqnEntranceTimeDef}
\sigma_t(x) = \inf\{ s | s \in[0,t] \mbox{ and } \forall r \in
(s,t], A_{r,t}(x) \in D \}.
\end{equation}
Let $g\dvtx\del D \times[0,\infty) \to\R$ and $\theta_0\dvtx D \to\R$ be
two given (regular enough) functions and define the process $\theta_t$ by
%
\begin{equation}\label{eqnThetaBddDef}
\theta_t(x) =
\cases{
g_{\sigma_t(x)} \circ A_{\sigma_t(x), t} (x),& \quad if $\sigma
_t(x) > 0 $,\cr
\theta_0 \circ A_{0,t}(x), & \quad if $\sigma_t(x) = 0$.
}
\end{equation}
Note that when $\sigma_t(x) > 0$, equation \eqref{eqnThetaBddDef} is
consistent with \eqref{eqnThetaDefAst}. Thus, \eqref{eqnThetaBddDef}
is the natural generalization of \eqref{eqnThetaDef} in the presence
of spatial boundaries and we expect $\bar\theta_t = E \theta_t$
satisfies the PDE \eqref{eqnBarTheta} with initial data $\bar\theta
_0 = \theta_0$ and boundary conditions $\theta= g$ on $\del D \times
[0, \infty)$. Indeed, this is essentially the expected value
representation obtained via the Kolmogorov backward equations.

If an extra term $c_t(x) \bar\theta_t(x)$ is desired on the left-hand
side of \eqref{eqnBarTheta}, then we only need to replace \eqref
{eqnThetaBddDef} by
\[
\theta_t(x) =
\cases{
\displaystyle \exp\biggl(-\int_{\sigma_t(x)}^t c_s(A_{s,t}) \,ds \biggr) g_{\sigma_t(x)} \circ
A_{\sigma_t(x), t} (x), &\quad  if $\sigma_t(x) > 0$, \cr
\displaystyle \exp\biggl(-\int_{0}^t c_s(A_{s,t}) \,ds \biggr) \theta_0 \circ A_{0,t}(x), & \quad
if $\sigma_t(x) = 0$
}
\]
provided $c$ is bounded below. This is essentially the Feynman--Kac
formula and its application to the Navier--Stokes equations is
developed in the next section.

\sidenote{Perhaps we can give meaning to \eqref{eqnTheta} if we only
differentiate after the backward exit time? It would be very
interesting to have a similar equation in the case with boundaries.}%
Note that the backward exit time $\sigma$ is usually discontinuous in
the spatial variable. Thus, even with smooth $g, \theta_0$, the
process $\theta$ need not be spatially continuous. As mentioned
earlier, equation \eqref{eqnTheta} will now become meaningless and we
will not be able to obtain a SPDE for $\theta$. However, equation~%
\eqref{eqnBarTheta}, which describes the evolution of the expected
value $\bar\theta= E \theta$, can be directly derived using the
backward Markov property and It\^o's formula (see, e.g.,~\cite
{bblFriedman}). We will not provide this proof here but will instead
provide a~proof for the more complicated analogue for the
Navier--Stokes equations described subsequently.

\subsection{Application to the Navier--Stokes equations in domains
with boundary}\label{sxnNSIntro}

First note that if $g = 0$ in \eqref{eqnThetaBddDef}, then the
solution to \eqref{eqnBarTheta} with initial data $\theta_0$ and
$0$-Dirichlet boundary conditions will be given by
%
\begin{equation}\label{eqnBarThetaBddDef}
\bar\theta_t = E \Chi{\sigma_t = 0} \theta_0 \circ A_{0,t}
\qquad \bigl[\mbox{i.e., }
\bar\theta_t(x) = E \Chi{\sigma_t(x) = 0} \theta_0 \circ
A_{0,t}(x) \bigr].
\end{equation}

Recall the no-slip boundary condition for the Navier--Stokes equations
is exactly a $0$-Dirichlet boundary condition on the velocity field.
Let $u$ be a~solution to the Navier--Stokes equations in $D$ with
initial data $u_0$ and no-slip boundary conditions. Now, following
\eqref{eqnBarThetaBddDef}, we would expect that analogous to \eqref
{eqnAveragedWebber}, the velocity field $u$ can be recovered from the
flow $A_{s,t}$ [equation~\eqref{eqnAst}], the backward exit time
$\sigma_t$ [equation~\eqref{eqnEntranceTimeDef}] and the initial data
$u_0$ by
%
\begin{equation}\label{eqnUNoBoundarysGen}
u_t = \lhp E \Chi{\sigma_t = 0} (\gradt A_{0,t}) u_0 \circ A_{0,t}.
\end{equation}

This, however, is \textit{false}. In fact, there are two elementary
reasons one should expect \eqref{eqnUNoBoundarysGen} to be false.
First, absorbing Brownian motion at the boundaries will certainly
violate incompressibility. The second and more fundamental reason is
that experiments and physical considerations lead us to expect
production of vorticity at the boundary. This is exactly what is
missing from \eqref{eqnUNoBoundarysGen}. The correct representation is
provided in the following result.

\begin{theorem}\label{thmSLNSNoSlip}
Let $u \in C^1([0,T); C^2(D)) \cap C([0,T]; C^1(\bar D))$ be a solution
of the Navier--Stokes equations \eqref{eqnNavierStokes} and \eqref{eqnIncompressibility} with initial data $u_0$ and no-slip boundary
conditions. Let $A$ be the solution to the backward SDE \eqref{eqnAst}
and $\sigma$ be the backward exit time defined by \eqref
{eqnEntranceTimeDef}. There exists a function $\tilde w\dvtx\del D \times
[0,T] \to\R^3$ such that for
%
\begin{equation}\label{eqnWdef}
w_t(x) =
\cases{
(\gradt A_{0,t}(x)) u_0 \circ A_{0,t}(x), & \quad when $\sigma_t =
0$,\cr
\bigl(\gradt A_{\sigma_t(x),t}(x)\bigr) \tilde w_{\sigma_t(x)} \circ A_{\sigma
_t(x),t}(x), &\quad  when $\sigma_t > 0$,
}
\end{equation}
we have
%
\begin{equation}\label{eqnUNoBoundarysCorrect}
u_t = \lhp E w_t.
\end{equation}

Conversely, given a function $\tilde w\dvtx \del D \times[0, T] \to\R
^d$, suppose there exists a~solution to the stochastic
system \eqref{eqnAst}, \eqref{eqnWdef}, \eqref
{eqnUNoBoundarysCorrect}. If further $u \in C^1([0,T);\break C^2(D)) \cap
C([0,T]; C^1(\bar D))$, then $u$ satisfies the Navier--Stokes
equa-\break tions~\mbox{\eqref{eqnNavierStokes}--\eqref{eqnIncompressibility}} with initial data $u_0$ and \textit
{vorticity} boundary conditions
%
\begin{equation}\label{eqnUVorticityBC}
\curl u = \curl E w \qquad \mbox{on }\del D \times[0, T].
\end{equation}
\end{theorem}

The proof of Theorem \ref{thmSLNSNoSlip} is presented in Section \ref
{sxnNoSlipProof}. We conclude this section with a few remarks.
\begin{remark}\label{rmkDiffNotation}
$\!\!\!\!\!$By $\!\gradt\! A_{\sigma_t(x), t}(x)$ in equation \eqref{eqnWdef} we mean
$[\!\gradt\! A_{s,t}(x)]_{s = \sigma_t(x)}$. That is, $\gradt A_{\sigma
_t(x), t}(x)$ refers to the transpose of the Jacobian of $A$, evaluated
at initial time $\sigma_t(x)$, final time $t$ and position $x$ (see
\cite
{bblKrylovQuasiDerivatives1,bblKrylovQuasiDerivatives2,bblKunita} for
existence). This is different from the transpose of the Jacobian of the
function~$A_{\sigma_t(\cdot),t}(\cdot)$ which does not exist as the
function is certainly not differentiable in space.
\end{remark}

\begin{remark}[(Regularity assumptions)]\label{rmkRegularityAssumptions}
In order to simplify the presentation, our regularity assumptions on
$u$ are somewhat generous. Our assumptions on $u$ will immediately
guarantee that $u$ has a Lipschitz extension to $\R^d$.
Now the process~$A$, defined to be a solution to \eqref{eqnAst} with
this Lipschitz extension of $u$, can be chosen to be a (backward)
stochastic flow of diffeomorphisms \cite{bblKunita}. Thus, $\grad A$
is well defined and further defining~$\sigma$ by~\eqref
{eqnEntranceTimeDef} is valid. Finally, since the statement of Theorem
\ref{thmSLNSNoSlip} only uses values of $A_{s,t}$ for $s \geq\sigma
_t$, the choice of the Lipschitz extension of $u$ will not matter. See
also Remark \ref{rmkWeakerRegularity}.
\end{remark}

\begin{remark}\label{rmkVorticityBC}
Note that our statement of the converse above \textit{does not}
explicitly give any information on the Dirichlet boundary values of
$u$. Of course, the normal component of $u$ must vanish at the boundary
of $D$ since~$u$ is the Leray--Hodge projection of a function. But an
explicit local relation between $\tilde w$ and the boundary values of
the tangential component of $u$ cannot be established. We remark,
however, that while the vorticity boundary condition \eqref
{eqnUVorticityBC} is somewhat artificial, it is enough to guarantee
uniqueness of solutions to the initial value problem for the
Navier--Stokes equations.
\end{remark}

\begin{remark}[(Choice of $\tilde w$)]\label{rmkTildeWChoice}
We explain how $\tilde w$ can be chosen to obtain the no-slip boundary
conditions. We will show (Lemma \ref{lmaWEquation}) that for $w$
defined by \eqref{eqnWdef}, the expected value $\bar w \defeq E w$
solves the PDE\sidenote{TODO: Cite some literature on Gage equations}
%
\begin{equation}\label{eqnBarWPDE}
\del_t \bar w_t + (u_t \cdot\grad) \bar w_t - \nu\lap\bar w_t +
(\gradt u_t) \bar w_t = 0
\end{equation}
with initial data
%
\begin{equation}\label{eqnBarWInitialData}
\bar w |_{t=0} = u_0.
\end{equation}
As shown before, $\gradt u_t$ in \eqref{eqnBarWPDE} denotes the
transpose of the Jacobian of $u_t$. Now, if $u = \lhp\bar w$, then we
will have $\curl u = \curl\bar w$ in $D$ and by continuity, on the
boundary of $D$. Thus, to prove existence of the function $\tilde w$,
we solve the PDE \eqref{eqnBarWPDE} with initial conditions \eqref
{eqnBarWInitialData} and \textit{vorticity} boundary conditions
%
\begin{equation}\label{eqnBarWBoundaryConditions}
\curl\bar w_t = \curl u_t \qquad \mbox{on }\del D.
\end{equation}
We chose $\tilde w$ to be the Dirichlet boundary values of this solution.
\end{remark}

To elaborate on Remark \ref{rmkTildeWChoice}, we trace through the
influence of the vorticity on the boundary on the velocity in the
interior. First, the vorticity at the boundary influences $\bar w$ by
entering as a boundary condition on the first derivatives for the PDE
\eqref{eqnBarWPDE}. Now, to obtain $u$ we need to find $\tilde w$, the
(Dirichlet) boundary values of \eqref{eqnBarWPDE} and use this to
weight trajectories that start on the boundary of $D$. The process of
finding $\tilde w$ is essentially passing from Neumann boundary values
of a PDE to the Dirichlet boundary values which is usually a nonlocal
pseudo-differential operator. Thus, while the procedure above is
explicit enough, the boundary vorticity influences the interior
velocity in a highly implicit, nonlocal manner.

\begin{remark}[(Uniqueness of $\tilde w$)]
Our choice of $\tilde w$ is not unique. Indeed, if $\bar w^1$ and $\bar
w^2$ are two solutions of \eqref{eqnBarWPDE}--\eqref
{eqnBarWBoundaryConditions}, then we must have $\bar w^1 - \bar w^2 =
\grad q$, where $q$ satisfies the equation
%
\begin{equation}\label{eqnQeq}
\grad\bigl( \del_t q + (u \cdot\grad) q - \nu\lap q \bigr) = 0
\end{equation}
with initial data $\grad q_0 = 0$. Since we do not have boundary
conditions on $q$, we can certainly have nontrivial solutions to this
equation. Thus, our choice of $\tilde w$ is only unique up to addition
by the gradient of a solution to \eqref{eqnQeq}.
\end{remark}
%

\section{Backward It\^o integrals}\label{sxnBackwardIto}
While the formulation of Theorem \ref{thmSLNSNoSlip} involves only
regular (forward) It\^o integrals, the proof requires backward It\^o
integrals and processes adapted to a two parameter filtration. The need
for backward It\^o integrals stems from equation \eqref{eqnAst} which,
as mentioned earlier, is the evolution of $A$, \textit{backward} in
time. This is, however, obscured because our diffusion coefficient is
constant making the martingale term exactly the increment of the Wiener
process and can be explicitly computed without any backward (or even
forward) It\^o integrals.

To elucidate matters, consider the flow $X'$ given by
%
\begin{equation}\label{eqnXstSigma}
X'_{s,t}(a) = a + \int_s^t u_r \circ X'_{s,r}(a) \,dr + \int_s^t \sigma
_r \circ X'_{s,r}(a) \,dW_r.
\end{equation}
If, as usual, $A'_{s,t} = (X'_{s,t})\inv$, then substituting formally%
\footnote{The formal substitution does \textit{not} give the correct
answer when $\sigma$ is not spatially constant. This is explained
subsequently and the correct equation is \eqref{eqnAstSigmaCorrect} below.}
$a = A'_{s,t}(x)$ and assuming the semigroup property gives the equation
%
\begin{equation}\label{eqnAstSigma}
A'_{s,t}(x) = x - \int_s^t u_r \circ A'_{r,t} (x) \,dr - \int_s^t
\sigma_r \circ A'_{r,t}(x) \,dW_r
\end{equation}
for the process $A'_{s,t}$. The need for backward It\^o integrals is
now evident; the last term above does not make sense as a forward It\^o
integral since $A'_{r,t}$ is not $\mathcal F_r$ measurable. This term,
however, is well defined as a backward It\^o integral; an integral with
respect to a decreasing filtration where processes are sampled at the
right endpoint. Since forward It\^o integrals are more predominant in
the literature, we recollect a few standard facts about backward It\^o
integrals in this section. A~more detailed account, with proofs, can be
found in \cite{bblFriedman,bblKunita}, for instance.

Let $(\Omega, \mathcal F, P)$ be a probability space, $\{W_t\}_{t\geq
0}$ be a $d$-dimensional Wiener process on $\Omega$ and let $\mathcal
F_{s,t}$ be the $\sigma$-algebra generated by the increments $W_{t'} -
W_{s'}$ for all $s\leq s' \leq t' \leq t$, augmented so that the
filtration $\{\mathcal F_{s,t}\}_{0\leq s \leq t}$ satisfies the usual
conditions\rlap{.}%
\footnote{By ``usual conditions'' in this context, we mean that for
all $s \geq0$, $\mathcal F_{s,s}$ contains all $\mathcal F_{0,\infty
}$-null sets. Further, $\mathcal F_{s,t}$ is right-continuous in $t$
and left-continuous in $s$. See \cite{bblKaratzasShreve}, {Definition
2.25,} for instance.}
Note that for $s \leq s' \leq t' \leq t$, we have $\mathcal
F_{s', t'} \subset\mathcal F_{s,t}$. Also $W_t - W_s$
is $\mathcal F_{s,t}$-measurable and is independent of both the past~%
$\mathcal F_{0,s}$, and the future $\mathcal F_{t,\infty}$.

We define a (two parameter) family of random variables $\{\xi_{s,t}\}
_{0 \leq s \leq t }$ to be a (two parameter) process adapted to the
(two parameter) filtration $\{\mathcal F_{s,t}\}_{0\leq s\leq t}$, if
for all $0 \leq s \leq t$, the random variable $\xi_{s,t}$ is
$\mathcal F_{s,t}$-measurable. For example, $\xi_{s,t} = W_t - W_s$ is
an adapted process. More generally, if $u$ and~$\sigma$ are regular
enough \textit{deterministic} functions, then the solution $\{
X'_{s,t}\}_{0 \leq s \leq t}$ of the (forward) SDE \eqref{eqnXstSigma}
is an adapted process.

Given an adapted (two parameter) process $\xi$ and any $t \geq0$, we
define the backward It\^o integral $\int_\cdot^t \xi_{r,t} \,dW_r$ by
\[
\int_s^t \xi_{r,t} \,dW_r = \lim_{\Vert P\Vert \to0} \sum_i \xi
_{t_{i+1}, t} (W_{t_{i+1}} - W_{t_i}),
\]
where $P = (r = t_0 < t_1 \cdots< t_N = t)$ is a partition of $[r, t]$
and $\Vert P\Vert $ is the length of the largest subinterval of $P$. The
limit is taken in the $L^2$ sense, exactly as with forward It\^o
integrals (see, e.g., {\cite{bblKaratzasShreve}, {page 148},
\cite{bblMcKean}, {page 35}, \cite{bblKunita}, {page 111}}).

The standard properties (existence, It\^o isometry, martingale
properties) of the backward It\^o integral are, of course, identical to
those of the forward integral. The only difference is in the sign of
the It\^o correction. Explicitly, consider the process $\{A'_{s,t}\}_{0
\leq s \leq t}$ satisfying the backward It\^o differential equation
\eqref{eqnAstSigma}. If $\{f_{s,t}\}_{0\leq s \leq t}$ is adapted,
$C^2$ in space and continuously differentiable with respect to $s$,
then the process $B_{s,t} = f_{s,t}\circ A_{s,t}$ satisfies the
backward It\^o differential equation
\begin{eqnarray*}
B_{t,t} - B_{s,t} &=& \int_s^t \biggl[\del_r f_{r,t} + (u_r \cdot\grad)
f_{r,t} - \frac{1}{2} a_r^{ij} \del_{ij} f_{r,t} \biggr]\circ{A_{r,t}} \,dr
\\
&&{}+\int_s^t [ \grad f_{r,t} \sigma_r ]\circ{A_{r,t}} \,dW_r,
\end{eqnarray*}
where $a_r^{ij} = \sigma_r^{ik} \sigma_r^{jk}$ with the Einstein sum
convention.

Though we only consider solutions to \eqref{eqnXstSigma} for constant
diffusion coefficient, we briefly address one issue when $\sigma$ is
not constant. Our motivation for the equation \eqref{eqnAstSigma} was
to make the substitution $x = A'_{s,t}(x)$ and formally use the
semigroup property. This, however, does not yield the correct equation
when $\sigma$ is not constant and the equation for $A'_{s,t} =
(X'_{s,t})\inv$ involves an additional correction term. To see this,
we discretize the forward integral in \eqref{eqnXstSigma} (in time)
and substitute $a = A'_{s,t}(x)$. This yields a sum sampled at the
\textit{left} endpoint of each time step. While this causes no
difficulty for the bounded variation terms, the martingale term is a
discrete approximation to a \textit{backward} integral and hence, must
be sampled at the \textit{right} endpoint of each time step. Converting
this to sum sampled at the right endpoint via a~Taylor expansion of
$\sigma$ is what gives this extra correction. Carrying through this
computation (see, e.g., \cite{bblKunita}, {Section 4.2}) yields the equation
\begin{eqnarray} \label{eqnAstSigmaCorrect}
A'_{s,t}(x) &=& x - \int_s^t u_r \circ A'_{r,t} (x) \,dr - \int_s^t
\sigma_r \circ A'_{r,t}(x) \,dW_r  \nonumber\\[-8pt]\\[-8pt]
&&{}+ \int_s^t \bigl( \del_j \sigma_r^{i,k} \circ A'_{r,t}(x) \bigr) \bigl( \sigma
^{j,k}_r \circ A'_{r,t}(x) \bigr) e_i \,dr,\nonumber
\end{eqnarray}
where $\{e_i\}_{1 \leq i \leq d}$ are the elementary basis vectors and
$\sigma^{i,j}$ denotes the $i, j$th entry in the $d \times d$
matrix $\sigma$.

We recall that the proof of the (forward) It\^o formula involves
approximating $f$ by its Taylor polynomial about the left endpoint of
the partition intervals. Analogously, the backward It\^o formula
involves approximating $f$ by Taylor polynomial about the \textit
{right} endpoint of partition intervals, which accounts for the
reversed sign in the It\^o correction.

Finally, we remark that for any fixed $t \geq0$, the solution $\{
A_{s,t}\}_{0\leq s \leq t}$ of the backward SDE \eqref{eqnAst} is a
backward strong Markov process [the same is true for solutions to
\eqref{eqnAstSigmaCorrect}]. The backward Markov property states that
$r < s < t$ then
\[
E_{\mathcal{F}_{s,t}} f\circ A_{r,t}(x) = E_{A_{s,t}(x)} f\circ A_{r,t}(x) = [ E f\circ A_{r,s}(y) ]_{y = A_{s,t}(x)},
\]
where $E_{\mathcal F_{s,t}}$ denotes the conditional expectation with
respect to the $\sigma$-alge\-bra $\mathcal F_{s,t}$ and
$E_{A_{s,t}(x)}$ the conditional expectation with respect to the
$\sigma$-alge\-bra generated by the process $A_{s,t}(x)$.

For the strong Markov property (we define $\sigma$ to be a \textit{backward $t$-stopping time}%
\footnote{Our use of the term backward $t$-stopping time is analogous
to $s$-stopping time in \cite{bblFriedman}, {page 24}.}
if almost surely $\sigma\leq t)$ and for all $s \leq t$, the event $\{
\sigma\geq s \}$ is $\mathcal F_{s,t}$-measurable. Now if $\sigma$ is
any backward $t$-stopping time with $r \leq\sigma\leq t$ almost
surely, the backward strong Markov property states
\[
E_{\mathcal F_{\sigma,t}} f\circ A_{r,t}(x) = E_{A_{\sigma,t}} f\circ A_{r,t}(x) = [ E f\circ A_{r,s}(y) ]\mathop{{}_{\hspace*{-21pt}
s = \sigma, }}_{y = A_{\sigma,t}(x).
}
\]
The proofs of the backward Markov properties is analogous to the proof
of the forward Markov properties and we refer the reader to \cite
{bblFriedman}, for instance.


\section{The no-slip boundary condition}\label{sxnNoSlipProof}
In this section we prove Theo\-rem~\ref{thmSLNSNoSlip}. First, we know
from \cite{bblKrylovQuasiDerivatives1,bblKrylovQuasiDerivatives2}
that spatial derivatives of $A$ can be interpreted as the limit (in
probability) of the usual difference quotient. In fact, for regular
enough velocity fields $u$ (extended to all of $\R^d$), the process
$A$ can, in fact, be chosen to be a flow of diffeomorphisms of $\R^d$
(see, e.g., \cite{bblKunita}) in which case $A$ is surely
differentiable in space. Interpreting the Jacobian of~$A$ as either the
limit (in probability) of the usual difference quotient or as the
Jacobian of the stochastic flow of diffeomorphism, we know \cite
{bblKrylovQuasiDerivatives1,bblKrylovQuasiDerivatives2,bblKunita} that~%
$\grad A$ satisfies the equation
%
\begin{equation}\label{eqnGradAst}
\grad A_{s,t}(x) = I - \int_s^t \grad u_r |_{A_{r,t}(x)} \grad
A_{r,t}(x)\, dr,
\end{equation}
obtained by formally differentiating \eqref{eqnAst} in space. Here $I$
denotes the $d \times d$ identity matrix. We reiterate that equation
\eqref{eqnGradAst} is an ODE as the Wiener process is independent of
the spatial parameter.


\begin{lemma}\label{lmaWEquation}
Let $D, u, T$ be as in Theorem \ref{thmSLNSNoSlip}, $\sigma$ be the
backward exit time from $D$ [equation \eqref{eqnEntranceTimeDef}]
and $A$ be the solution to \eqref{eqnAst} with respect to the backward
stopping time $\sigma$.
\begin{enumerate}[(2)]
\item[(1)] Let $\bar w \in C^1([0,T); C^2(D)) \cap C([0,T]; C^1(\bar D))$ be
the solution of \eqref{eqnBarWPDE} with initial data \eqref
{eqnBarWInitialData} and boundary conditions
%
\begin{equation}\label{eqnBarWDirichletBC}
\bar w = \tilde w \qquad \mbox{on } \del D.
\end{equation}
Then, for $w$ defined by \eqref{eqnWdef}, we have $\bar w = E w$.

\item[(2)] Let $w$ be defined by \eqref{eqnWdef} and $\bar w = E w$ as
above. If for all $t \in(0, T]$, $\bar w_t \in\mathcal D(A_{\cdot,
t})$ and $\bar w$ is $C^1$ in time, then $\bar w$ satisfies
%
\begin{equation}\label{eqnBarWPDEwithL}
\del_t \bar w + L_t \bar w + (\gradt u) \bar w = 0,
\end{equation}
where $L_t$ is defined by
%
\begin{equation}\label{eqnLdef}
L_t \phi(x) = \lim_{s \to t^-} \frac{\phi(x) - E \phi( A_{s\varmax
\sigma_t(x),t}(x) )}{t-s}
\end{equation}
and $\mathcal D(A_{\cdot, t})$ is the set of all $\phi$ for which the
limit on the right-hand side exists. Further, $\bar w$ has initial data
$u_0$ and boundary conditions \eqref{eqnBarWDirichletBC}.
\end{enumerate}
\end{lemma}

Before proceeding any further, we first address the relationship
between the two assertions of the lemma. We claim that if $\bar w \in
C^1((0,T); C^2(D))$, then equation~\eqref{eqnBarWPDEwithL} reduces to
equation \eqref{eqnBarWPDE}. This follows immediately from the next
proposition.
\begin{proposition}\label{ppnGeneratorAt}
If $\phi\in C^2(D)$, then for any $t \in(0, T]$, $\phi\in\mathcal
D(A_{\cdot, t})$ and further, $L_t\phi= (u_t \cdot\grad) \phi- \nu
\lap\phi$.
\end{proposition}
%

\begin{pf}
Omitting the spatial variable for notational convenience, the backward
It\^o formula gives
\begin{eqnarray*}
\phi- \phi\circ A_{s\varmax\sigma_t, t} &=& \phi\circ A_{t,t} -
\phi\circ A_{s\varmax\sigma_t, t}\\
&=& \int_{s\varmax\sigma_t}^t [ (u_r \cdot\grad) \phi|
_{A_{r,t}} - \nu\lap\phi\at_{A_{r,t}} ] \,dr + \sqrt{2\nu} \int
_{s \varmax\sigma_t}^t \grad\phi\at_{A_{r,t}} \,dW_r.
\end{eqnarray*}
Since $s \varmax\sigma_t$ is a backward $t$-stopping time, the second
term above is a~martingale. Thus
\begin{eqnarray*}
L_t\phi&=& \lim_{s \to t^-} E \frac{1}{t-s} \int_s^t \Chi{r \geq
\sigma_t} [ (u_r \cdot\grad) \phi\at_{A_{r,t}} - \nu\lap\phi
\at_{A_{r,t}} ] \,dr\\
&=& (u_t \cdot\grad)\phi- \nu\lap\phi
\end{eqnarray*}
since the process $A$ has continuous paths and $\sigma_t < t$ on the
interior of $D$.
\end{pf}
%

\begin{remark}\label{rmkWeakerRegularity}
One can weaken the regularity assumptions on $u$ in the statement of
Theorem \ref{thmSLNSNoSlip} by instead assuming for all $t \in(0,
T]$, $u_t \in\mathcal D(A_{\cdot, t})$ and is $C^1$ in time, as with
the second assertion of Lemma \ref{lmaWEquation}. However, while the
formal calculus remains essentially unchanged, there are a couple of
technical points that require attention. First, when assumptions on
smoothness of $u$ up to the boundary is relaxed (or when $\del D$ is
irregular), a~Lipschitz extension of $u$ need not exist. In this case,
we can no longer use~\eqref{eqnEntranceTimeDef} to define $\sigma$.
Further, we can not regard the process $A$ as a stochastic flow of
diffeomorphisms and some care has to be taken when differentiating it.
These issues can be addressed using relatively standard techniques and
once they are sorted out, the proof of Theorem~\ref{thmSLNSNoSlip}
remains unchanged.
\end{remark}

Now we prove the first assertion of Lemma \ref{lmaWEquation}.

\begin{pf}
$\!\!$Recall that $\gradt A_{s,t}$ is differentiable in $s$. Differentiating
\eqref{eqnGradAst} in~$s$ and transposing the matrices gives
%
\begin{equation}\label{eqnGradtAst}
\del_s \grad A_{s,t}(x) = \gradt A_{s,t}(x) \gradt u_s \at_{A_{s,t}(x)}.
\end{equation}
Let $t \in(0, T]$, $x \in D$ and $\sigma'$ be any backward
$t$-stopping time with $\sigma' \geq\sigma_t(x)$ almost surely.
Omitting the spatial variable for convenience, the backward It\^o
formula and equations \eqref{eqnBarWPDE} and \eqref{eqnGradtAst} give
\begin{eqnarray*}
&&\bar w_t - \gradt A_{\sigma', t} \bar w_{\sigma'}
\circ A_{\sigma', t} \\
&&\qquad = \gradt A_{t,t} \bar w_t \circ A_{t,t} - \gradt A_{\sigma
', t} \bar w_{\sigma'} \circ A_{\sigma', t}\\
&&\qquad = \int_{\sigma'}^t \del_r \gradt A_{r,t} \bar w_r \circ A_{r,t}
\\
&&\qquad \quad {} + \int_{\sigma'}^t \gradt A_{r,t} \bigl( \del_r \bar w_r + (u_r \cdot
\grad) \bar w_r - \nu\lap\bar w_r \bigr) \circ A_{r,t}\, dr  \\
&&\quad \qquad {} + \sqrt{2\nu} \int_{\sigma'}^t (\gradt A_{r,t}) (\gradt\bar w_r)
\circ A_{r,t} \,dW_r\\
&&\qquad = \int_{\sigma'}^t \gradt A_{r,t} \bigl( (\gradt u_r) \bar w_r + \del_r
\bar w_r + (u_r \cdot\grad) \bar w_r - \nu\lap\bar w_r \bigr) \circ
A_{r,t} \,dr  \\
&&\quad\qquad {}  + \sqrt{2\nu} \int_{\sigma'}^t (\gradt A_{r,t}) (\gradt\bar w_r)
\circ A_{r,t} \,dW_r\\
&&\qquad = \sqrt{2\nu} \int_{\sigma'}^t (\gradt A_{r,t}) w_r \circ A_{r,t} \,dW_r.
\end{eqnarray*}
Thus, taking expected values gives
%
\begin{equation}\label{eqnBarWSigma}
\bar w_t(x) = E \gradt A_{\sigma', t}(x) \bar w_{\sigma'} \circ A_{\sigma', t}(x).
\end{equation}
Recall that when $\sigma_t(x) > 0$, $A_{\sigma_t(x), t}(x) \in\del
D$. Thus, choosing $\sigma' = \sigma_t(x)$ and using the boundary
conditions \eqref{eqnBarWDirichletBC} and initial data \eqref
{eqnBarWInitialData}, we have
%
\begin{equation}\label{eqnWSigma}
\bar w_{\sigma_t(x)} \circ A_{\sigma_t(x), t} =
\cases{
\tilde w_{\sigma_t(x)} \circ A_{\sigma_t(x), t}, &\quad  if $\sigma
_t(x) > 0$,\cr
u_0 \circ A_{\sigma_t(x), t}, &\quad  if $\sigma_t(x) = 0$.
}
\end{equation}
Substituting this in \eqref{eqnBarWSigma} completes the proof.
\end{pf}

In order to prove the second assertion in Lemma \ref{lmaWEquation}, we
will directly prove \eqref{eqnBarWSigma} using the backward strong
Markov property. Before beginning the proof, we establish a few preliminaries.

Let $D$, $u$, $T$, $\sigma$, $A$, $w$, $\bar w$ be as in the second
assertion of Lemma \ref{lmaWEquation}. Given $x \in D$ and a $d \times
d$ matrix $M$, define the process $\{B_{s,t}(x,M)\}_{\sigma_t(x) \leq
s \leq t \leq T}$ to be the solution of the ODE
\[
B_{s,t}(x,M) = M - \int_s^t \grad u_r \at_{A_{r,t}(x)} B_{r,t}(x,M) \,dr.
\]
If $I$ denotes the $d \times d$ identity matrix, then by \eqref
{eqnGradAst} we have $B_{s,t}(x,I) = \grad A_{s,t}(x)$ for any $\sigma
_t(x) \leq s \leq t \leq T$. Further, since the evolution equation for~%
$B$ is linear, we see
%
\begin{equation}\label{eqnBstM}
B_{s,t}(x, M) = B_{s,t}(x,I) M = \grad A_{s,t}(x) M.
\end{equation}
Note that for any fixed $t \in(0, T]$, the process $\{\grad A_{s,t}\}
_{0 \leq s \leq t}$ is \textit{not} a backward Markov process. Indeed,
the evolution of $\grad A_{s,t}$ at any time $s \leq t$ depends on the
time $s$ through the process $A_{s,t}$ appearing on the right-hand side
in~\eqref{eqnGradAst}. However, process $(A_{s,t}, \grad A_{s,t})$ [or
equivalently the process $(A_{s,t}, B_{s,t})$] is a backward Markov
process since the evolution of this system now only depends on the
state. This leads us to the following identity which is the essence of
proof of the second assertion in Lemma \ref{lmaWEquation}.\vspace*{-3pt}

\begin{lemma}\label{lmaMarkov}
Choose any backward $t$-stopping time $\sigma'$ with $\sigma' \geq
\sigma_t(x)$ almost surely. Then
\begin{eqnarray} \label{eqnMarkov}
&&E_{\mathcal F_{\sigma', t}} B_{\sigma_t(x),
t}^*(x, I) \bar w_{\sigma_t(x)} \circ A_{\sigma_t(x),
t}(x)\nonumber\\[-8pt]\\[-8pt]
&&\qquad = \mathop{\bigl[ E B_{\sigma_r(y), r}^*(y, M) \bar w_{\sigma_r(y)} \circ A_{\sigma
_r(y), r}(y) \bigr]_{
r = \sigma', y = A_{\sigma
', t}(x),}}_{\hspace*{150pt} M = B_{\sigma', t}(x,I)}\nonumber
\end{eqnarray}
holds almost surely.\vspace*{-3pt}
\end{lemma}

This follows from an appropriate application of the backward strong
Markov property. While this is easily believed, checking that the
strong Markov property applies in this situation requires a little work
and will distract from the heart of the matter. Thus, we momentarily
postpone the proof of Lemma \ref{lmaMarkov} and proceed with the proof
of the second assertion of Lemma \ref{lmaWEquation}.\vspace*{-3pt}

\begin{pf*}{Proof of Lemma \ref{lmaWEquation}}
We recall $\bar w = Ew$ where $w$ is defined by \eqref{eqnWdef}. By
our assumption on $u$ and $\del D$, the boundary conditions \eqref
{eqnBarWDirichletBC} and initial data \eqref{eqnBarWInitialData} are
satisfied. For convenience, when $y \in\del D$, $t > 0$, we define
$w_t(y) = \tilde w(y)$ and when $t = 0$, $y \in\bar D$, we define
$w_0(y) = u_0(y)$.

Let $x \in D$, $t \in (0, T]$ as used before. Let $\sigma'$ be any
backward $t$-stopping time with $\sigma' \geq\sigma_t(x)$. First, if
$\sigma' = \sigma_t(x)$ almost surely, then, since the point
$(A_{\sigma_t(x), t}, t)$ belongs to the parabolic boundary $\del_p
(D \times[0, T]) \defeq(\del D \times[0, T]) \cup(D \times\{0\})$,
our boundary conditions and initial data will guarantee~\eqref{eqnBarWSigma}.

Now, for arbitrary $\sigma' \geq\sigma_t(x)$, we will use Lemma \ref
{lmaMarkov} to deduce \eqref{eqnBarWSigma} directly. Indeed,
\begin{eqnarray*}
\bar w_t(x) &=& E \gradt A_{\sigma_t(x), t}(x) \bar w_{\sigma_t(x)}
\circ A_{\sigma_t(x), t}\\
&=& E E_{\mathcal F_{\sigma', t}} B_{\sigma
_t(x), t}^*(x, I) \bar w_{\sigma_t(x)} \circ A_{\sigma_t(x), t}(x)\\
&=& E \bigl( \mathop{\bigl[ E B_{\sigma_r(y), r}^*(y, M) \bar w_{\sigma_r(y)} \circ
A_{\sigma_r(y), r}(y) \bigr]_{
r = \sigma', y = A_{\sigma
', t}(x),}}_{\hspace*{150pt} M = B_{\sigma', t}(x,I)
} \bigr)\\
&=& E \bigl( \mathop{\bigl[ M^* E B_{\sigma_r(y), r}^*(y, I) \bar w_{\sigma_r(y)} \circ
A_{\sigma_r(y), r}(y) \bigr]_{
r = \sigma', y = A_{\sigma', t}(x),}}_{\hspace*{160pt} M = B_{\sigma', t}(x,I)
} \bigr)\\
&=& E \gradt A_{\sigma',t}(x) \bar w_{\sigma'} \circ
A_{\sigma', t}(x),
\end{eqnarray*}
showing that \eqref{eqnBarWSigma} holds for any backward $t$ stopping
time $\sigma' \geq\sigma_t(x)$.

Now, choose $\sigma' = s \varmax\sigma_t(x)$ for $s < t$. Note that
for any $x \in D$, we must have $\sigma_t(x) < t$ almost surely. Thus,
omitting the spatial coordinate for convenience, we have
\begin{eqnarray*}
0 &=& \lim_{s \to t^-} \frac{ \bar w_t - \bar w_t }{t - s} = \lim_{s
\to t^-} \frac{1}{t-s} ( \bar w_t - E \gradt A_{s \varmax\sigma_t,
t} \bar w_{s \varmax\sigma_t} \circ A_{s \varmax\sigma_t, t} ) \\
&= &\lim_{s \to t^-} \biggl(
\frac{1}{t-s} [ \bar w_t - E \bar w_t \circ A_{s \varmax\sigma_t,
t} ]  \\
&&\hphantom{\lim_{s \to t^-} \biggl(}
{} + \frac{1}{t-s} E (\bar w_t - \bar w_{s \varmax\sigma_t} ) \circ
A_{s \varmax\sigma_t, t} \\
&&\hphantom{\lim_{s \to t^-} \biggl(}
{} + \frac{1}{t-s} E ( I - \gradt A_{s \varmax\sigma_t, t} ) \bar
w_{s \varmax\sigma_t} \circ A_{s \varmax\sigma_t, t} \biggr)
\\
&=& L_t \bar w_t + \del_t \bar w_t + (\gradt u_t) \bar w_t,
\end{eqnarray*}
on the interior of $D$. The proof is complete.
\end{pf*}
%

It remains to prove Lemma \ref{lmaMarkov}.
$\!\!\!$\begin{pf*}{Proof of Lemma \ref{lmaMarkov}}
$\!\!\!$Define the stopped processes $A'_{s,t}(x) \!=\! A_{\sigma_t(x) \varmax s,
t}(x)$ and $B'_{s,t}(x, M) = B_{\sigma_t(x) \varmax s, t}(x, M)$.
Define the process $C$ by
\[
C_{s,t}(x, M, \tau) = \bigl(A'_{s,t}(x), B'_{s, t}( x, M), \tau+ t -
\sigma_t(x) \varmax s\bigr).
\]
Note that for any given $s \leq t$, we know that $\sigma_t(x)$ need
not be $\mathcal F_{s,t}$ measurable. However, $\sigma_t(x) \varmax s$
\textit{is} an $\mathcal F_{s,t}$ measurable backward $t$-stopping
time. Thus, $A'_{s,t}$, $B'_{s,t}$ and, consequently, $C_{s,t}$ are all
$\mathcal F_{s,t}$ measurable.

Now we claim that almost surely, for $0 \leq r \leq s \leq t \leq T$,
we have the backward semigroup identity
%
\begin{equation}\label{eqnCSemigroup}
C_{r,s} \circ C_{s,t} = C_{r,t}.
\end{equation}
To prove this, consider first the third component of the left-hand side
of~\eqref{eqnCSemigroup}:
%
\begin{equation}\label{eqnC3LHS}
C_{r,s}^{(3)} \circ C_{s,t}(x, M, \tau) = \bigl(\tau+ t - \sigma_t(x)
\varmax s\bigr) + s - \sigma_s(A'_{s,t}(x)) \varmax s.
\end{equation}
Consider the event $\{s > \sigma_t(x)\}$. By the semigroup property
for $A$ and strong existence and uniqueness of solutions to \eqref
{eqnAst}, we have $\sigma_s(A_{s,t}(x)) = \sigma_t(x)$ almost surely.
Thus, almost surely on $\{s > \sigma_t(x)\}$, we have
\begin{eqnarray*}
C_{r,s}^{(3)} \circ C_{s,t}(x, M, \tau) &=& (\tau+ t - s) + s - \sigma
_t(x) \varmax s\\
&=& \tau+ t - \sigma_t(x) \varmax r = C^{(3)}_{r,t}(x, M, \tau).
\end{eqnarray*}
Now consider the event $\{s \!\leq\!\sigma_t\}$. We know $A'_{s,t}(x)
\!\in\!\del D$ and so \mbox{$\sigma_s(A'_{s,t}(x)) \!=\! s$}. This gives
\[
C_{r,s}^{(3)} \circ C_{s,t}(x, M, \tau) = \bigl(\tau+ t - \sigma_t(x) \bigr) +
s - s = \tau+ t - \sigma_t(x) \varmax r = C^{(3)}_{r,t}(x)
\]
almost surely on $\{s \leq\sigma_t(x)\}$. Therefore, we have proved
almost sure equality of the third components in equation \eqref{eqnCSemigroup}.

For the first component $C^{(1)}_{s,t} = A'_{s,t}$, consider as before the case
$s > \sigma_t(x)$. In this case $A'_{s,t} = A_{s,t}$ and the semigroup
property of $A$ gives equality of the first components in \eqref
{eqnCSemigroup} almost surely on $\{s > \sigma_t(x)\}$. When $ s \leq
\sigma_t(x)$, as before, $A'_{s,t} \in\del D$ and $\sigma
_s(A'_{s,t}(x)) = s$. Thus,
\[
A'_{r,s} \circ A'_{s,t}(x) = A_{s,s} \circ A_{\sigma_t(x), t}(x) =
A_{\sigma_t(x), t}(x) = A'_{r,t}(x)
\]
almost surely on $s \leq\sigma_t(x)$. This shows almost sure equality
of the first components in equation \eqref{eqnCSemigroup}. Almost sure
equality of the second components follows similarly, completing the
proof of \eqref{eqnCSemigroup}.

Now, for $0 \leq r \leq s \leq t \leq T$, the random variable $C_{s,t}$
is $\mathcal F_{s,t}$ measurable and so must be independent of
$\mathcal F_{r,s}$. This, along with \eqref{eqnCSemigroup}, will
immediately guarantee the Markov property for $C$. Since the filtration
$\mathcal F_{\cdot, \cdot}$ satisfies the usual conditions and for
any fixed $t$ the function $s \mapsto C_{s,t}$ is continuous, $C$
satisfies the strong Markov property (see, e.g., \cite
{bblFriedman}, {Theorem 2.4}).

Thus, for any fixed $t \in[0, T]$ and any Borel function $\varphi$,
the strong Markov property gives
\begin{eqnarray*}
E_{\mathcal F_{\sigma', t} } \varphi( C_{0,t}(x, I,
0)) &=& [ E \varphi( C_{r,t}(y, M, \tau)) ]\mathop{{}_{\hspace*{-58pt}
r = \sigma',}}_{\hspace*{-1pt} (y, M, \tau) = C_{0,
\sigma'}(x, I, 0)
}\\
&=& [ E \varphi( C_{r,t}(y, M, \tau)) ]\mathop{{}_{\hspace*{-22pt}
r = \sigma', y = A_{\sigma
',t}(x),}}_{\hspace*{-1pt}M = B_{\sigma',t}(x, I), \tau
= \sigma_r(x),
}
\end{eqnarray*}
{\spaceskip=0.2em plus 0.05em minus 0.02em almost surely for any
$x \!\in\!\R^d$, $M \!\in\!\R^{d^2}$, $\tau\!\geq\!0$. Choosing $\varphi(x, M, \tau)
\!=\!M\transpose\bar w_{t - \tau}(x)$} proves \eqref{eqnMarkov}.
\end{pf*}
%

Now a direct computation shows that if $\bar w$ satisfies \eqref
{eqnBarWPDE}, then $u = \lhp\bar w$ satisfies \eqref{eqnNavierStokes}
regardless of our choice of $\tilde w$. Of course, we will only get the
no-slip boundary conditions with the correct choice of $\tilde w$. We
first obtain the PDE for $u$.
\begin{lemma}\label{lmaUPDE}
If $\bar w$ satisfies \eqref{eqnBarWPDE} and $u = \lhp\bar w$, then
$u$ satisfies \eqref{eqnNavierStokes} and~\eqref{eqnIncompressibility}.
\end{lemma}
%

\begin{pf}
By definition of the Leray--Hodge projection, $u = w + \grad q$ for
some function $q$ and equation \eqref{eqnIncompressibility} is
automatically satisfied. Thus, using equation \eqref{eqnBarWPDE} we have
\begin{eqnarray} \label{eqnuNS1}
&&\del_t u_t + (u_t \cdot\grad) u_t - \nu\lap u_t + (\gradt u_t) u_t
\nonumber\\[-8pt]\\[-8pt]
&&\qquad {}+ \del_t \grad q_t + (u_t \cdot\grad) \grad q_t + (\gradt u_t) \grad
q_t - \nu\lap\grad q_t = 0.\nonumber
\end{eqnarray}
Defining $p$ by
\[
\grad p = \grad\bigl( \tfrac{1}{2} |u|^2 + \del_t q_t + (u_t \cdot
\grad) q_t - \nu\lap q_t \bigr),
\]
equation \eqref{eqnuNS1} becomes \eqref{eqnNavierStokes}.
\end{pf}
%

Now to address the no-slip boundary condition. The curl of $\bar w$
satisfies the vorticity equation which is how the vorticity enters our
boundary condition.\vspace*{-3pt}
\begin{lemma}\label{lmaCurlW}
Let $\bar w$ be a solution of \eqref{eqnBarWPDE}. Then $\xi= \curl
\bar w$ satisfies the vorticity equation
%
\begin{equation}\label{eqnCurlW}
\del_t \xi+ (u \cdot\grad) \xi- \nu\lap\xi=
\cases{
0, &\quad  if $d = 2$,\cr
(\xi\cdot\grad) u, &\quad  if $d = 3$.
}
\end{equation}
\end{lemma}

\begin{pf}
We only provide the proof for $d=3$. For this proof we will use
subscripts to indicate the component instead of time as we usually do.
If $i,j,k \in\{1,2,3\}$ are all distinct, let $\epsilon_{ijk}$ denote
the signature of the permutation $(1, 2, 3) \mapsto(i,j,k)$. For
convenience, we let $\epsilon_{ijk} = 0$ if $i, j, k$ are not all
distinct. Using the Einstein summation convention, $\xi= \curl\bar w$
translates to $\xi_i = \epsilon_{ijk} \del_j \bar w_k$ on
components. Thus, taking the curl of \eqref{eqnBarWPDE} gives
%
\begin{equation}\label{eqnXi}
\del_t \xi_i + (u \cdot\grad) \xi_i - \nu\lap\xi_i + \epsilon
_{ijk} \,\del_j u_m \,\del_m \bar w_k + \epsilon_{ijk} \,\del_k u_m \,\del
_j \bar w_m = 0
\end{equation}
because $\epsilon_{ijk} \,\del_j \del_k u_m \bar w_m = 0$. Making the
substitutions $j \mapsto k$ and $k \mapsto j$ in the last sum above we have
\begin{eqnarray*}
\epsilon_{ijk} \,\del_j u_m \,\del_m \bar w_k + \epsilon_{ijk} \,\del_k
u_m \,\del_j \bar w_m &= &\epsilon_{ijk}\, \del_j u_m ( \del_m \bar w_k
- \del_k \bar w_m )\\
&= &\epsilon_{ijk} \,\del_j u_m \epsilon_{nmk} \xi_n\\
&=& ( \delta_{in} \delta_{jm} - \delta_{im} \delta_{jn} ) \,\del_j
u_m \xi_n\\
&=& - \del_j u_i \xi_j,
\end{eqnarray*}
where $\delta_{ij}$ denotes the Kronecker delta function and the last
equality follows because $\del_j u_j = 0$. Thus, \eqref{eqnXi}
reduces to \eqref{eqnCurlW}.\vspace*{-3pt}
\end{pf}

Theorem \ref{thmSLNSNoSlip} now follows from the above lemmas.
\begin{pf*}{Proof of Theorem \ref{thmSLNSNoSlip}}
First, suppose $u$ is a solution of the Navier--Stokes equations, as in
the statement of the theorem. We choose $\tilde w$ as explained in
Remark~\ref{rmkTildeWChoice}. Notice that our assumptions on $u$ and
$D$ will guarantee a classical solution to \eqref{eqnBarWPDE}--\eqref
{eqnBarWBoundaryConditions} exists on the interval $[0, T]$ and thus,
such a choice is possible.

By Lemma \ref{lmaWEquation} we see that for $w$ defined by \eqref
{eqnWdef}, the expected value $\bar w = E w$ satisfies \eqref
{eqnBarWPDE} with initial data \eqref{eqnBarWInitialData} and boundary
conditions~\eqref{eqnBarWDirichletBC}. By our choice of $\tilde w$ and
uniqueness to the Dirichlet problem \eqref{eqnBarWPDE}, \eqref
{eqnBarWInitialData} and \eqref{eqnBarWDirichletBC}, we must have the
vorticity boundary condition \eqref{eqnBarWBoundaryConditions}.

Now, let $\xi= \curl\bar w$ and $\omega= \curl u$. By Lemma \ref
{lmaCurlW}, we see that $\xi$ satisfies the vorticity equation \eqref
{eqnCurlW}. Since $u$ satisfies \eqref{eqnNavierStokes} and \eqref{eqnIncompressibility}, it is well known (see, e.g., \cite
{bblChorinMarsden,bblMajdaBertozzi} or the proof of Lemma \ref
{lmaCurlW}) that $\omega$ also satisfies
%
\begin{equation}\label{eqnVorticity}
\del_t \omega_t + (u_t \cdot\grad) \omega_t - \nu\lap\omega_t =
\cases{
0, &\quad  if $d = 2$,\cr
(\omega_t \cdot\grad) u_t, &\quad  if $d = 3$.
}
\end{equation}
From \eqref{eqnBarWBoundaryConditions} we know $\xi= \omega$ on
$\del D \times[0, T]$. By \eqref{eqnBarWInitialData}, we see that
$\xi_0 = \curl u_0 = \omega_0$ and hence, $\xi= \omega$ on the
parabolic boundary $\del_p (D \times[0, T])$.

The above shows that $\omega$ and $\xi$ both satisfy the same PDE
[equations~\eqref{eqnCurlW} or \eqref{eqnVorticity}] with the same
initial data and boundary conditions and so we must have $\xi= \omega
$ on $D \times[0, T]$. Thus, $\curl\bar w = \curl u$ in $D \times[0,
T]$ showing $u$ and~$\bar w$ differ by a gradient. Since $\divergence u
= 0$ and $u = 0$ on $\del D \times[0, T]$, we must have $u = \lhp\bar
w$ proving \eqref{eqnUNoBoundarysCorrect}.

Conversely, assume we have a solution to the system \eqref{eqnAst},
\eqref{eqnWdef} and~\eqref{eqnUNoBoundarysCorrect}. As stated above,
Lemma \ref{lmaWEquation} shows $\bar w = Ew$ satisfies~\eqref
{eqnBarWPDE} with initial data~\eqref{eqnBarWInitialData}. By Lemma
\ref{lmaUPDE} we know $u$ satisfies the equation \eqref{eqnNavierStokes} and \eqref{eqnIncompressibility}
with initial data
$u_0$. Finally, since equation \eqref{eqnUNoBoundarysCorrect} shows
$\curl u = \curl\bar w$ in $D \times[0, T]$ and by continuity, we
have the boundary condition \eqref{eqnUVorticityBC}.
\end{pf*}
%

\section{Vorticity transport and ideally conserved quantities}\label{sxnVorticityTransport}
The vorticity is a quantity which is of fundamental importance, both
for the physical and theoretical aspects of fluid dynamics. To single
out one among the numerous applications of vorticity, we refer the
reader to two classical criterion which guarantee global and existence
and regularity of the Navier--Stokes equations provided the vorticity
is appropriately controlled: the first due to Beale, Kato and Majda
\cite{bblBealeKatoMajda} and the second due to Constantin and
Fefferman \cite{bblConstFefferman}.

For the Euler equations, exact identities and conservation laws
governing the evolution of vorticity are well known. For instance,
vorticity transport [equation \eqref{eqnEulerVortTransport}] shows
that the vorticity at time $t$ followed along streamlines is exactly
the initial vorticity stretched by the Jacobian of the flow map.
Similarly, the conservation of circulation [equation \eqref
{eqnEulerCirculation}] shows that the line integral of the velocity
(which, by Stokes theorem, is a surface integral of the vorticity)
computed along a closed curve that is transported by the fluid flow is
constant in time.

Prior to \cite{bblConstIyerSLNS}, these identities were unavailable
for the Navier--Stokes equations. In \cite{bblConstIyerSLNS}, the
authors provide analogues of these identities for the Navier--Stokes
equations in the absence of boundaries. These identities, however, do
not always prevail in the presence of boundaries.

In this section we illustrate the issues involved by considering three
inviscid identities. All three identities generalize perfectly to the
viscous situations without boundaries. In the presence of boundaries,
the first identity (vorticity transport) generalizes perfectly, the
second identity (Ertel's Theorem) generalizes somewhat unsatisfactorily
and the third identity (conservation of circulation) has no nontrivial
generalization in the presence of boundaries.

\subsection{Vorticity transport}
Let $u^0$ be a solution to the Euler equations with initial data $u_0$.
Let $X^0$ the inviscid flow map defined by \eqref{eqnEulerXdef} and
for any $t \geq0$, let $A^0_t = (X^0_t)\inv$ be the spatial inverse
of the diffeomorphism $X^0_t$. The vorticity transport (or Cauchy
formula) states
%
\begin{equation}\label{eqnEulerVortTransport}
\omega^0_t =
\cases{
\omega^0_0 \circ A^0_t, &\quad  if $d = 2$,\cr
[ (\grad X^0_t) \omega^0_0 ] \circ A^0_t, &\quad  if $d = 3$,
}
%
\end{equation}
where we recall that the vorticity $\omega^0$ is defined by $\omega^0
= \curl u^0$ and where $\omega^0_0 = \curl u_0$ is the initial vorticity.

In \cite{bblConstIyerSLNS}, the authors obtained a natural
generalization of \eqref{eqnEulerVortTransport} for the Navier--Stokes
equations in the absence of spatial boundaries. If $u$ solves~\eqref{eqnNavierStokes} and \eqref{eqnIncompressibility} with initial data $u_0$
and $X$ is the noisy flow map defied by \eqref{eqnXdef}--\eqref
{eqnX0}, then $\omega= \curl u$ is given by
%
\begin{equation}\label{eqnNSVortTransport}
\omega_t =
\cases{
E \omega_0 \circ A_t, &\quad  if $d = 2$,\cr
E ( (\grad X_t) \omega_0 ) \circ A_t, &\quad  if $d= 3$.
}
\end{equation}

We now provide the generalization of this in the presence of
boundaries. Note that for any $t \geq0$, $(\grad X_t) \circ A_t =
(\grad A_t)\inv$, so we can rewrite \eqref{eqnNSVortTransport}
completely in terms of the process $A$. Now, as usual, we replace $A =
X\inv$ with the solution of \eqref{eqnAst} with respect to the
minimal existence time $\sigma$. We recall that in Theorem \ref
{thmSLNSNoSlip}, in addition to ``starting trajectories at the
boundary,'' we had to correct the expression for the velocity
by the boundary values of a related quantity (the vorticity). For the
vorticity, however, we need no additional correction and the interior
vorticity is completely determined given~$A$, $\sigma$ and the
vorticity on the parabolic boundary%
\footnote{Recall the parabolic boundary $\del_p (D \times[0, T])$
is defined to be $(D \times\{0\}) \cup(\del D \times[0, T))$.}
$\del_p (D \times[0, T])$.

\begin{proposition}\label{ppnVorticityRep}
Let $u$ be a solution to \eqref{eqnNavierStokes} and \eqref{eqnIncompressibility} in $D$ with initial data $u_0$ and suppose
$\omega= \curl u \in C^1([0,T); C^2(D)) \cap C([0, T] \times\bar D)$.
Let $\tilde\omega$ denote the values of $\omega$ on the parabolic
boundary $\del_p (D \times[0, T])$. Explicitly, $\tilde\omega$ is
defined by
\[
\tilde\omega(x, t) =
\cases{
\omega_0(x), & \quad if $x \in D $ and $ t = 0$,\cr
\omega_t(x), & \quad if $x \in\del D$.
}
\]
Then,
%
\begin{equation}\label{eqnOmegaRep}
\omega_t(x) =
\cases{
E \bigl[ \tilde\omega_{\sigma_t(x)} \bigl( A_{\sigma_t(x), t}(x) \bigr) \bigr], &\quad
if $d = 2$,\cr
E \bigl[ \bigl(\grad A_{\sigma_t(x), t}(x) \bigr)\inv\tilde\omega_{\sigma_t(x)}
\bigl( A_{\sigma_t(x), t}(x) \bigr) \bigr], & \quad if $ d = 3$.
}
\end{equation}
\end{proposition}

\begin{proposition}\label{ppnVorticityEq}
More generally, suppose $\tilde\omega$ is any function defined on the
parabolic boundary of $D \times[0, T]$ and let $\omega$ be defined by
\eqref{eqnOmegaRep}. If for all $t \in(0, T]$, $\omega_t \in
\mathcal D(A_{\cdot, t})$ and $\omega$ is $C^1$ in time, then $\omega
$ satisfies
\[
\del_t \omega_t + L_t \omega_t =
\cases{
0, &\quad  if $d = 2$,\cr
(\omega_t \cdot\grad) u_t, &\quad  if $d = 3$,
}
\]
with $\omega= \tilde\omega$ on the parabolic boundary. Here, $L_t$
is the generator of $A_{\cdot, t}$; $\mathcal D(A_{\cdot, t})$ is the
domain of $L_t$. These are defined in the statement of Lemma~\ref
{lmaWEquation}.
\end{proposition}
%

Of course, Proposition \ref{ppnVorticityEq}, along with Proposition
\ref{ppnGeneratorAt} and uniqueness of (strong) solutions to \eqref
{eqnVorticity}, will prove Proposition \ref{ppnVorticityRep}. However,
direct proofs of both Proposition \ref{ppnVorticityEq} and Proposition
\ref{ppnVorticityRep} are short and instructive and we provide
independent proofs of each.

\begin{pf*}{Proof of Proposition \ref{ppnVorticityRep}}
We only provide the proof when $d = 3$. As shown before,
differentiating \eqref{eqnGradAst} in space and taking the matrix
inverse of both sides gives
%
\begin{equation}\label{eqnDsGradAstInv}
\del_r (\grad A_{r,t}(x) )\inv= - (\grad A_{r,t}(x) )\inv\grad
u_r\at_{A_{r,t}(x)},
\end{equation}
almost surely. Now choose any $x \in D$, $t > 0$ and any backward
$t$-stopping time $\sigma' \geq\sigma_t(x)$. Omitting the spatial
parameter for notational convenience, the backward It\^o formula gives
\begin{eqnarray*}
&&\omega_t - (\grad A_{\sigma',t} )\inv
\omega_{\sigma'} \circ A_{\sigma',t} \\
&&\qquad  = (\grad A_{t,t} )\inv\omega_t \circ A_{t,t} - (\grad A_{\sigma
',t} )\inv\omega_{\sigma'} \circ A_{\sigma
',t}\\
&&\qquad =
\int_{\sigma'}^t \del_r (\grad A_{r,t} )\inv\omega_r
\circ A_{r,t} \,dr  \\
&&\qquad \quad
{}+ \int_{\sigma'}^t (\grad A_{r,t} )\inv\bigl(\del_r \omega_r + (u_r
\cdot\grad) \omega_r - \nu\lap\omega_r \bigr) \circ A_{r,t} \,dr  \\
&&\qquad \quad
{}+ \sqrt{2\nu} \int_{\sigma'}^t (\grad A_{r,t} )\inv(\grad\omega
_r ) \circ A_{r,t} \,dW_r
\\
&&\quad =
\int_{\sigma'}^t - (\grad A_{r,t} )\inv\grad u_r\at
_{A_{r,t}} \omega_r \circ A_{r,t} \,dr \\
&&\quad \qquad {} + \int_{\sigma'}^t (\grad A_{r,t} )\inv((\omega_r \cdot\grad)
u_r ) \circ A_{r,t} \,dr \\
&&\quad \qquad {}  + \sqrt{2\nu}\int_{\sigma'}^t (\grad A_{r,t} )\inv(\grad\omega
_r ) \circ A_{r,t}\, dW_r
\\
&&\qquad =\sqrt{2\nu}\int_{\sigma'}^t (\grad A_{r,t} )\inv(\grad\omega
_r ) \circ A_{r,t} \,dW_r.
\end{eqnarray*}
Thus, taking expected values gives
%
\begin{equation}\label{eqnOmegaMarkov}
\omega_t = E [ (\grad A_{\sigma',t})\inv\omega
_{\sigma'} \circ A_{\sigma',t} ].
\end{equation}
Choosing $\sigma' = \sigma_t(x)$ and using the fact that $A_{\sigma
_t(x),t}(x)$ always belongs to the parabolic boundary finishes the proof.
\end{pf*}

\begin{pf*}{Proof of Proposition \ref{ppnVorticityEq}}
Again, we only consider the case \mbox{$d = 3$}. We will prove equation \eqref
{eqnOmegaMarkov} directly and then deduce \eqref{eqnVorticity}. Let
the process~$B$ be as in the proof of the second assertion of Lemma
\ref{lmaWEquation} and use $B\inv$ to denote the process consisting
of \textit{matrix} inverses of the process $B$. Pick $x \in D$, $t \in
(0, T]$ and a backward $t$-stopping time $\sigma' \geq\sigma_t(x)$.
Using \eqref{eqnMarkov} we have
\begin{eqnarray*}
\omega_t(x) &=& E \bigl[ \bigl(\grad A_{\sigma_t(x), t}(x) \bigr)\inv\tilde\omega
_{\sigma_t(x)} \bigl( A_{\sigma_t(x), t}(x) \bigr) \bigr]\\
&=& E E_{\mathcal F_{\sigma',t}} \bigl[ B_{\sigma_t(x),t}\inv
(x,I) \tilde\omega_{\sigma_t(x)} \circ A_{\sigma_t(x),t}(x) \bigr]\\
&=& E \bigl( \bigl[ E B_{\sigma_r(y),r}\inv(y,M) \tilde\omega_{\sigma_r(y)}
\circ A_{\sigma_r(y),r}(y) \bigr]\mathop{{}_{
r = \sigma', y = A_{\sigma',t}(x),}}_{\hspace*{-15pt}M=B_{\sigma',t}(x,I)
} \bigr)\\
&=& E \bigl( \bigl[ M\inv E B_{\sigma_r(y),r}\inv(y,I) \tilde\omega_{\sigma
_r(y)} \circ A_{\sigma_r(y),r}(y) \bigr]\mathop{{}_{
r = \sigma', y = A_{\sigma',t}(x),}}_{\hspace*{-15pt}M=B_{\sigma',t}(x,I)
} \bigr)\\
&=& E [ (\grad A_{\sigma',t}(x) )\inv\omega_{\sigma'} \circ
A_{\sigma',t}(x) ],
\end{eqnarray*}
proving \eqref{eqnOmegaMarkov}.

As stated before, choose $s \leq t$ and $\sigma' = \sigma_t(x)
\varmax s$. Omitting the spatial parameter for notational convenience gives
\begin{eqnarray*}
0 &=& \lim_{s \to t^-} \frac{\omega_t - \omega_t}{t-s} = \lim_{s
\to t^-} \frac{1}{t-s} [ \omega_t - E (\grad A_{\sigma_t \varmax
s,t} )\inv\omega_{\sigma_t\varmax s} \circ A_{\sigma_t\varmax s,t}
]\\
&= &\lim_{s \to t^-} \biggl(
\frac{1}{t-s} [ \omega_t - E\omega_t \circ A_{\sigma_t
\varmax s, t} ] \\
&&\hphantom{\lim_{s \to t^-} \biggl(}
{} + \frac{1}{t-s} E [ \omega_t - \omega_{\sigma_t\varmax s} ] \circ
A_{\sigma_t\varmax s,t} \\
&&\hphantom{\lim_{s \to t^-} \biggl(}
{} + \frac{1}{t-s} E [ I - (\grad A_{\sigma_t \varmax s,t} )\inv]
\omega_{\sigma_t\varmax s} \circ A_{\sigma_t\varmax s,t} \biggr)
\\
&=&L_t \omega_t + \del_t \omega_t - (\grad u_t) \omega_t.
\end{eqnarray*}
\upqed
\end{pf*}

%

We remark that the vorticity transport in Propositions \ref
{ppnVorticityRep} or \ref{ppnVorticityEq} can be used to provide a
stochastic representation of the Navier--Stokes equations. To see this,
first note that the proofs of Propositions \ref{ppnVorticityRep} and
\ref{ppnVorticityEq} are independent of Theorem \ref{thmSLNSNoSlip}.
Next, since $u$ is divergence free, taking the curl twice gives the
negative laplacian. Thus, provided boundary conditions on~$u$ are
specified, we can obtain $u$ from $\omega$ by
%
\begin{equation}\label{eqnBiotSavart}
u_t = (- \lap)\inv\curl\omega_t.
\end{equation}
Therefore, in Theorem \ref{thmSLNSNoSlip} we can replace \eqref
{eqnUNoBoundarysCorrect} by \eqref{eqnOmegaRep} and \eqref
{eqnBiotSavart}, where~$\tilde\omega$ is the vorticity on the
parabolic boundary and we impose $0$-Dirichlet boundary conditions on
\eqref{eqnBiotSavart}.

\subsection{Ertel's theorem}

As shown above, we use a superscript of $0$ to denote the appropriate
quantities related to the Euler equations. For this section we also
assume $d = 3$. Ertel's theorem says that if $\theta^0$ is constant
along trajectories of $X^0$, then so is $(\omega^0 \cdot\grad)
\theta^0$. Hence, $\phi^0 = (\omega^0 \cdot\grad) \theta^0$
satisfies the PDE
\[
\del_t \phi^0 + (u \cdot\grad) \phi^0 = 0.
\]

For the Navier--Stokes equations, we first consider the situation
without boundaries. Let $u$ solve \eqref{eqnNavierStokes} and \eqref{eqnIncompressibility}, $X$ be defined by \eqref{eqnXdef}, $A$ be the
spatial inverse of $X$ and define $\xi$ by
\[
\xi_t(x) = (\grad A_t(x))\inv\omega_0 \circ A_t(x),
\]
where $\omega_0 = \curl u_0$ is the initial vorticity. From \eqref
{eqnNSVortTransport} we know that $\omega= \curl u = E \xi$. Now we
can generalize Ertel's theorem as follows:

\begin{proposition}\label{ppnErtelNoBoundaries}
Let $\theta$ be a $C^1(\R^d)$ valued process. If $\theta$ is
constant along trajectories of the (stochastic) flow $X$, then so is
$(\xi\cdot\grad) \theta$. Hence, $\phi= E (\xi\cdot\grad)
\theta$ satisfies the PDE
%
\begin{equation}\label{eqnPhiPDE}
\del_t \phi_t + (u_t \cdot\grad) \phi_t - \nu\lap\phi_t = 0,
\end{equation}
with initial data $(\omega_0 \cdot\grad) \theta_0$.
\end{proposition}
\begin{pf}
If $\theta$ is constant along trajectories of $X$, we must have
$\theta_t = \theta_0 \circ A_t$ almost surely. Thus,
\[
(\xi_t \cdot\grad) \theta_t = (\grad\theta_t) \xi_t = \grad
\theta_0 \at_{A_t} (\grad A_t) (\grad A_t)\inv\omega_0 \circ A_t =
( \xi_0 \cdot\grad\theta_0 ) \circ A_t,
\]
which is certainly constant along trajectories of $X$. The PDE for
$\phi$ now follows immediately.
\end{pf}

Now, in the presence of boundaries, this needs further modification.
Let~$A$ be a solution to \eqref{eqnAst} and $\sigma$ be the backward
exit time of $A$ from $D$. The notion of ``constant along
trajectories'' now corresponds to processes $\theta$ defined~by
%
\begin{equation}\label{eqnThetaBddDef1}
\theta_t(x) = \tilde\theta_{\sigma_t(x)}\bigl(A_{\sigma_t(x), t}\bigr),
\end{equation}
for some function $\tilde\theta$ defined on the parabolic boundary of $D$.

Irrespective of the regularity of $D$ and $\tilde\theta$, the process
$\theta$ will not be continuous in space, let alone differentiable.
The problem arises because while $A$ is regular enough in the spatial
variable, the existence time $\sigma_t$ is not. Thus, we are forced to
avoid derivatives on $\sigma$ in the statement of the theorem, leading
to a somewhat unsatisfactory generalization.

\begin{proposition}\label{ppnErtelBoundaries}
Let $\tilde\theta$ be a $C^1$ function defined on the parabolic
boundary of $D \times[0, T]$ and let $\tilde\theta'$ be any $C^1$
extension of $\tilde\theta$, defined in a neighborhood of the
parabolic boundary of $D \times[0, T]$. If $\theta$ is defined by
\eqref{eqnThetaBddDef1}, then
\[
\phi_t(x) = E [(\xi_t \cdot\grad) (\tilde\theta'_s \circ
A_{s,t})(x) ]_{s = \sigma_t(x)}
\]
satisfies the PDE \eqref{eqnPhiPDE} with initial data $(\omega_0
\cdot\grad) \tilde\theta_0$ and boundary conditions $\phi_t(x) =
(\omega_t \cdot\grad) \tilde\theta'(x)$ for $x \in\del D$.
\end{proposition}

\begin{remark*}
A satisfactory generalization in the scenario with boundaries would be
to make sense of $E (\xi_t \cdot\grad) \theta_t$ (despite the
spatial discontinuity of~$\theta$) and reformulate Proposition \ref
{ppnErtelBoundaries} accordingly.
\end{remark*}

Note that when $D = \R^d$, then $\sigma_t \equiv0$ and
hence, $\phi_t = E (\xi_t \cdot\grad) \theta_t$. In this case
Proposition \ref{ppnErtelBoundaries} reduces to Proposition~\ref
{ppnErtelNoBoundaries}. The proof of Proposition~\ref
{ppnErtelBoundaries} is identical to that of Proposition~\ref
{ppnErtelNoBoundaries} and the same argument obtains
\[
[(\xi_t \cdot\grad) (\tilde\theta'_s \circ A_{s,t})(x) ]_{s =
\sigma_t(x)} = [ (\xi_s \cdot\grad) \tilde\theta'_s (y) ]\mathop{{}_{\hspace*{-16pt}
s = \sigma_t(x),}}_{y=A_{\sigma_t(x),t}(x)},
\]
which immediately implies \eqref{eqnPhiPDE}.

\subsection{Circulation}

The circulation is the line integral of the velocity field along a
closed curve. For the Euler equations, the circulation along a closed
curve that is transported by the flow is constant in time. Explicitly,
let $u^0$, $X^0$, $A^0$, $u_0$ be as in the previous subsection. Let
$\Gamma$ be a rectifiable closed curve, then for any $t \geq0$,
%
\begin{equation}\label{eqnEulerCirculation}
\oint_{X^0_t(\Gamma)} u^0_t \cdot dl = \oint_{\Gamma} u^0_0 \cdot dl.
\end{equation}
For the Navier--Stokes equations, without boundaries, a generalization
of~\eqref{eqnEulerCirculation} was considered in \cite
{bblConstIyerSLNS}. Let $u$ solve \eqref{eqnNavierStokes} and \eqref{eqnIncompressibility}, $X$ be defined by~\eqref{eqnXdef} and \eqref
{eqnX0} and $A$ be the spatial inverse of $X$. Then
%
\begin{equation}\label{eqnNSCirculation}
\oint_\Gamma u_t \cdot dl = E\oint_{A_t(\Gamma)} u_0 \cdot dl.
\end{equation}

A proof of this (in the absence of boundaries) follows immediately from
Theorem \ref{thmSLNS}. Indeed,
\begin{eqnarray} \label{eqnConsCirculationProof}
E\oint_{A_t(\Gamma)} u_0 \cdot dl &=& E\oint_{\Gamma} (\gradt A_t )
u_0 \circ A_t \cdot dl \nonumber\\[-8pt]\\[-8pt]
&=& E \oint_\Gamma\lhp[ (\gradt A_t ) u_0 \circ A_t ] \cdot dl =
\oint_\Gamma u_t \cdot dl,\nonumber
\end{eqnarray}
where the first equality follows by definition of line integrals, the
second because the line integral of gradients along closed curves is
$0$ and the last by Fubini and \eqref{eqnAveragedWebber}.

Equation \eqref{eqnNSCirculation} does not make sense in the presence
of boundaries, as the curves one integrates over will no longer be rectifiable!


\section*{Acknowledgment}
The authors would like to thank the referee for suggesting various
improvements to the first version of the paper.

%

\printaddresses


\begin{thebibliography}{40}

\bibitem{bblArnold}
%
\begin{barticle}[mr]
\bauthor{\bsnm{Arnold}, \bfnm{V.}\binits{V.}}
(\byear{1966}).
\btitle{Sur la g\'eom\'etrie diff\'erentielle des groupes de {L}ie de dimension
infinie et ses applications \`a l'hydrodynamique des fluides parfaits}.
\bjournal{Ann. Inst. Fourier (Grenoble)}
\bvolume{16}
\bpages{319--361}.
\bid{mr={0202082}}
\end{barticle}
%
\endbibitem

\bibitem{bblBealeKatoMajda}
%
\begin{barticle}[mr]
\bauthor{\bsnm{Beale}, \bfnm{J. T.}\binits{J. T.}},
\bauthor{\bsnm{Kato}, \bfnm{T.}\binits{T.}} \AND
\bauthor{\bsnm{Majda}, \bfnm{A.}\binits{A.}}
(\byear{1984}).
\btitle{Remarks on the breakdown of smooth solutions for the {$3$}-{D} {E}uler
equations}.
\bjournal{Comm. Math. Phys.}
\bvolume{94}
\bpages{61--66}.
\bid{mr={0763762}}
\end{barticle}
%
\endbibitem

\bibitem{bblBhattacharyaORST1}
%
\begin{barticle}[mr]
\bauthor{\bsnm{Bhattacharya}, \bfnm{Rabi N.}\binits{R. N.}},
\bauthor{\bsnm{Chen}, \bfnm{Larry}\binits{L.}},
\bauthor{\bsnm{Dobson}, \bfnm{Scott}\binits{S.}},
\bauthor{\bsnm{Guenther}, \bfnm{Ronald B.}\binits{R. B.}},
\bauthor{\bsnm{Orum}, \bfnm{Chris}\binits{C.}},
\bauthor{\bsnm{Ossiander}, \bfnm{Mina}\binits{M.}},
\bauthor{\bsnm{Thomann}, \bfnm{Enrique}\binits{E.}} \AND
\bauthor{\bsnm{Waymire}, \bfnm{Edward C.}\binits{E. C.}}
(\byear{2003}).
\btitle{Majorizing kernels and stochastic cascades with applications to
incompressible {N}avier--{S}tokes equations}.
\bjournal{Trans. Amer. Math. Soc.}
\bvolume{355}
\bpages{5003--5040 (electronic)}.
\bid{doi={10.1090/S0002-9947-03-03246-X}, mr={1997593}}
\end{barticle}
%
\endbibitem

\bibitem{bblBhattacharyaORST2}
%
\begin{bincollection}[mr]
\bauthor{\bsnm{Bhattacharya}, \bfnm{Rabi}\binits{R.}},
\bauthor{\bsnm{Chen}, \bfnm{Larry}\binits{L.}},
\bauthor{\bsnm{Guenther}, \bfnm{Ronald B.}\binits{R. B.}},
\bauthor{\bsnm{Orum}, \bfnm{Chris}\binits{C.}},
\bauthor{\bsnm{Ossiander}, \bfnm{Mina}\binits{M.}},
\bauthor{\bsnm{Thomann}, \bfnm{Enrique}\binits{E.}} \AND
\bauthor{\bsnm{Waymire}, \bfnm{Edward C.}\binits{E. C.}}
(\byear{2005}).
\btitle{Semi-{M}arkov cascade representations of local solutions to 3-{D}
incompressible {N}avier--{S}tokes}.
In \bbooktitle{Probability and Partial Differential Equations in
Modern Applied
Mathematics}
(\beditor{E. C. Waymire and J. Duan}, eds.).
\bseries{IMA Vol. Math. Appl.}
\bvolume{140}
\bpages{27--40}.
\bpublisher{Springer}, \baddress{New York}.
\bid{doi={10.1007/978-0-387-29371-4_3}, mr={2202031}}
\end{bincollection}
%
\endbibitem

\bibitem{bblBusnello}
%
\begin{barticle}[mr]
\bauthor{\bsnm{Busnello}, \bfnm{Barbara}\binits{B.}}
(\byear{1999}).
\btitle{A probabilistic approach to the two-dimensional {N}avier--{S}tokes
equations}.
\bjournal{Ann. Probab.}
\bvolume{27}
\bpages{1750--1780}.
\bid{doi={10.1214/aop/1022677547}, mr={1742887}}
\end{barticle}
%
\endbibitem

\bibitem{bblBusnelloFlandolliRomito}
%
\begin{barticle}[mr]
\bauthor{\bsnm{Busnello}, \bfnm{Barbara}\binits{B.}},
\bauthor{\bsnm{Flandoli}, \bfnm{Franco}\binits{F.}} \AND
\bauthor{\bsnm{Romito}, \bfnm{Marco}\binits{M.}}
(\byear{2005}).
\btitle{A probabilistic representation for the vorticity of a three-dimensional
viscous fluid and for general systems of parabolic equations}.
\bjournal{Proc. Edinb. Math. Soc. (2)}
\bvolume{48}
\bpages{295--336}.
\bid{doi={10.1017/S0013091503000506}, mr={2157249}}
\end{barticle}
%
\endbibitem

\bibitem{bblConstELE}
%
\begin{barticle}[mr]
\bauthor{\bsnm{Constantin}, \bfnm{Peter}\binits{P.}}
(\byear{2001}).
\btitle{An {E}ulerian-{L}agrangian approach for incompressible fluids: Local
theory}.
\bjournal{J. Amer. Math. Soc.}
\bvolume{14}
\bpages{263--278 (electronic)}.
\bid{doi={10.1090/S0894-0347-00-00364-7}, mr={1815212}}
\end{barticle}
%
\endbibitem

\bibitem{bblConstEntropies}
%
\begin{barticle}[mr]
\bauthor{\bsnm{Constantin}, \bfnm{Peter}\binits{P.}}
(\byear{2006}).
\btitle{Generalized relative entropies and stochastic representation}.
\bjournal{Int. Math. Res. Not.}
\bpages{Art. ID 39487, 9}.
\bid{doi={10.1155/IMRN/2006/39487}, mr={2250023}}
\end{barticle}
%
\endbibitem

\bibitem{bblConstFefferman}
%
\begin{barticle}[mr]
\bauthor{\bsnm{Constantin}, \bfnm{Peter}\binits{P.}} \AND
\bauthor{\bsnm{Fefferman}, \bfnm{Charles}\binits{C.}}
(\byear{1993}).
\btitle{Direction of vorticity and the problem of global regularity
for the
{N}avier--{S}tokes equations}.
\bjournal{Indiana Univ. Math. J.}
\bvolume{42}
\bpages{775--789}.
\bid{doi={10.1512/iumj.1993.42.42034}, mr={1254117}}
\end{barticle}
%
\endbibitem

\bibitem{bblConstFoias}
%
\begin{bbook}[mr]
\bauthor{\bsnm{Constantin}, \bfnm{Peter}\binits{P.}} \AND
\bauthor{\bsnm{Foias}, \bfnm{Ciprian}\binits{C.}}
(\byear{1988}).
\btitle{Navier--{S}tokes Equations}.
\bpublisher{Univ.  Chicago Press}, \baddress{Chicago, IL}.
\bid{mr={0972259}}
\end{bbook}
%
\endbibitem

\bibitem{bblConstIyerSLNS}
%
\begin{barticle}[mr]
\bauthor{\bsnm{Constantin}, \bfnm{Peter}\binits{P.}} \AND
\bauthor{\bsnm{Iyer}, \bfnm{Gautam}\binits{G.}}
(\byear{2008}).
\btitle{A stochastic {L}agrangian representation of the three-dimensional
incompressible {N}avier--{S}tokes equations}.
\bjournal{Comm. Pure Appl. Math.}
\bvolume{61}
\bpages{330--345}.
\bid{doi={10.1002/cpa.20192}, mr={2376844}}
\end{barticle}
%
\endbibitem

\bibitem{bblConstIyerSEntropy}
%
\begin{barticle}[mr]
\bauthor{\bsnm{Constantin}, \bfnm{Peter}\binits{P.}} \AND
\bauthor{\bsnm{Iyer}, \bfnm{Gautam}\binits{G.}}
(\byear{2006}).
\btitle{Stochastic {L}agrangian transport and generalized relative entropies}.
\bjournal{Commun. Math. Sci.}
\bvolume{4}
\bpages{767--777}.
\bid{mr={2264819}}
\end{barticle}
%
\endbibitem

\bibitem{bblCiprianoCruzerio}
%
\begin{barticle}[mr]
\bauthor{\bsnm{Cipriano}, \bfnm{F.}\binits{F.}} \AND
\bauthor{\bsnm{Cruzeiro}, \bfnm{A. B.}\binits{A. B.}}
(\byear{2007}).
\btitle{Navier-{S}tokes equation and diffusions on the group of homeomorphisms
of the torus}.
\bjournal{Comm. Math. Phys.}
\bvolume{275}
\bpages{255--269}.
\bid{doi={10.1007/s00220-007-0306-3}, mr={2335775}}
\end{barticle}
%
\endbibitem

\bibitem{bblChorin}
%
\begin{barticle}[mr]
\bauthor{\bsnm{Chorin}, \bfnm{Alexandre Joel}\binits{A. J.}}
(\byear{1973}).
\btitle{Numerical study of slightly viscous flow}.
\bjournal{J. Fluid Mech.}
\bvolume{57}
\bpages{785--796}.
\bid{mr={0395483}}
\end{barticle}
%
\endbibitem

\bibitem{bblChorinMarsden}
%
\begin{bbook}[mr]
\bauthor{\bsnm{Chorin}, \bfnm{Alexandre J.}\binits{A. J.}} \AND
\bauthor{\bsnm{Marsden}, \bfnm{Jerrold E.}\binits{J. E.}}
(\byear{1993}).
\btitle{A Mathematical Introduction to Fluid Mechanics},
\bedition{3rd} ed.
\bseries{Texts in Applied Mathematics}
\bvolume{4}.
\bpublisher{Springer}, \baddress{New York}.
\bid{mr={1218879}}
\end{bbook}
%
\endbibitem

\bibitem{bblEyink}
%
\begin{barticle}[mr]
\bauthor{\bsnm{Eyink}, \bfnm{Gregory L.}\binits{G. L.}}
(\byear{2010}).
\btitle{Stochastic least-action principle for the incompressible
{N}avier--{S}tokes equation}.
\bjournal{Phys. D}
\bvolume{239}
\bpages{1236--1240}.
\bid{doi={10.1016/j.physd.2008.11.011}, mr={2657460}}
\bptnote{check year}
\end{barticle}
%
\endbibitem

\bibitem{bblEyinkSL}
%
\begin{bmisc}[auto:SpringerTagBib|2010-03-24|17:41:21]
\bauthor{\bsnm{Eyink}, \bfnm{Gregory L.}\binits{G. L.}} (\byear{2008}).
\bhowpublished{{Stochastic line-motion and stochastic conservation
laws for
non-ideal hydromagnetic models. I. Incompressible fluids and isotropic
transport coefficients}. Preprint. Available at \href
{http://arxiv.org/abs/0812.0153}{arXiv:0812.0153}}.
\end{bmisc}
%
\endbibitem

\bibitem{bblFriedman}
%
\begin{bbook}[mr]
\bauthor{\bsnm{Friedman}, \bfnm{Avner}\binits{A.}}
(\byear{2006}).
\btitle{Stochastic Differential Equations and Applications}.
\bpublisher{Dover}, \baddress{Mineola, NY}.
(\bnote{Two volumes bound as one. Reprint of the 1975 and 1976 original
published in two volumes}.)
\bid{mr={2295424}}
\end{bbook}
%
\endbibitem

\bibitem{bblGliklikh}
%
\begin{bbook}[mr]
\bauthor{\bsnm{Gliklikh}, \bfnm{Yuri}\binits{Y.}}
(\byear{1997}).
\btitle{Global Analysis in Mathematical Physics: Geometric and Stochastic Methods}.
\bseries{Applied Mathematical Sciences}
\bvolume{122}.
\bpublisher{Springer}, \baddress{New York}.
(\bnote{Translated from the 1989 Russian
original and with Appendix F by Viktor L. Ginzburg}.)
\bid{mr={1438545}}
\end{bbook}
%
\endbibitem

\bibitem{bblIyerThesis}
%
\begin{bmisc}[auto:SpringerTagBib|2010-03-24|17:41:21]
\bauthor{\bsnm{Iyer}, \bfnm{Gautam}\binits{G.}} (\byear{2006}).
\bhowpublished{A stochastic Lagrangian formulation of the
Navier--Stokes and
related transport equations. Ph.D. thesis, Univ. Chicago}.
\end{bmisc}
%
\endbibitem

\bibitem{bblKaratzasShreve}
%
\begin{bbook}[mr]
\bauthor{\bsnm{Karatzas}, \bfnm{Ioannis}\binits{I.}} \AND
\bauthor{\bsnm{Shreve}, \bfnm{Steven E.}\binits{S. E.}}
(\byear{1991}).
\btitle{Brownian Motion and Stochastic Calculus},
\bedition{2nd} ed.
\bseries{Graduate Texts in Mathematics}
\bvolume{113}.
\bpublisher{Springer}, \baddress{New York}.
\bid{mr={1121940}}
\end{bbook}
%
\endbibitem

\bibitem{bblKrylovQuasiDerivatives1}
%
\begin{bincollection}[mr]
\bauthor{\bsnm{Krylov}, \bfnm{N. V.}\binits{N. V.}}
(\byear{1993}).
\btitle{Quasiderivatives for solutions of {I}t\^o's stochastic
equations and
their applications}.
In \bbooktitle{Stochastic Analysis and Related Topics ({O}slo, 1992)}.
\bseries{Stochastics Monogr.}
\bvolume{8}
\bpages{1--44}.
\bpublisher{Gordon and Breach}, \baddress{Montreux}.
\bid{mr={1268004}}
\end{bincollection}
%
\endbibitem

\bibitem{bblKrylovQuasiDerivatives2}
%
\begin{barticle}[vtex]
\bauthor{\bsnm{Krylov}, \bfnm{N. V.}\binits{N. V.}}
(\byear{2004}).
\btitle{Quasiderivatives and interior smoothness of harmonic functions
associated with degenerate diffusion processes}.
\bjournal{Electron. J. Probab.}
\bvolume{9}
\bpages{615--633 (electronic)}.
\bid{mr={2082053}}
\end{barticle}
%
\endbibitem


\bibitem{bblKuznetsovRuban}
%
\begin{barticle}[mr]
\bauthor{\bsnm{Kuznetsov}, \bfnm{E. A.}\binits{E. A.}} \AND
\bauthor{\bsnm{Ruban}, \bfnm{V. P.}\binits{V. P.}}
(\byear{2000}).
\btitle{Hamiltonian dynamics of vortex and magnetic lines in
hydrodynamic type
systems}.
\bjournal{Phys. Rev. E (3)}
\bvolume{61}
\bpages{831--841}.
\bid{doi={10.1103/PhysRevE.61.831}, mr={1736469}}
\end{barticle}
%
\endbibitem

\bibitem{bblKunita}
%
\begin{bbook}[mr]
\bauthor{\bsnm{Kunita}, \bfnm{Hiroshi}\binits{H.}}
(\byear{1997}).
\btitle{Stochastic Flows and Stochastic Differential Equations}.
\bseries{Cambridge Studies in Advanced Mathematics}
\bvolume{24}.
\bpublisher{Cambridge Univ. Press}, \baddress{Cambridge}.
(\bnote{Reprint of the 1990 original}.)
\bid{mr={1472487}}
\end{bbook}
%
\endbibitem

\bibitem{bblLeJanSznitman}
%
\begin{barticle}[mr]
\bauthor{\bsnm{Le Jan}, \bfnm{Y.}\binits{Y.}} \AND
\bauthor{\bsnm{Sznitman}, \bfnm{A. S.}\binits{A. S.}}
(\byear{1997}).
\btitle{Stochastic cascades and $3$-dimen\-sional {N}avier--{S}tokes
equations}.
\bjournal{Probab. Theory Related Fields}
\bvolume{109}
\mbox{\bpages{343--366}}.
\bid{doi={10.1007/s004400050135}, mr={1481125}}
\end{barticle}
%
\endbibitem

\bibitem{bblMcKean}
%
\begin{bbook}[mr]
\bauthor{\bsnm{McKean}, \bfnm{H. P.}\binits{H. P.} \bsuffix{Jr.}}
(\byear{1969}).
\btitle{Stochastic Integrals}.
\bseries{Probability and Mathematical Statistics}
\bvolume{5}.
\bpublisher{Academic Press}, \baddress{New York}.
\bid{mr={0247684}}
\end{bbook}
%
\endbibitem

\bibitem{bblMajdaBertozzi}
%
\begin{bbook}[mr]
\bauthor{\bsnm{Majda}, \bfnm{Andrew J.}\binits{A. J.}} \AND
\bauthor{\bsnm{Bertozzi}, \bfnm{Andrea L.}\binits{A. L.}}
(\byear{2002}).
\btitle{Vorticity and Incompressible Flow}.
\bseries{Cambridge Texts in Applied Mathematics}
\bvolume{27}.
\bpublisher{Cambridge Univ. Press}, \baddress{Cambridge}.
\bid{mr={1867882}}
\end{bbook}
%
\endbibitem

\bibitem{bblMichelMischlerPerthame1}
%
\begin{barticle}[mr]
\bauthor{\bsnm{Michel}, \bfnm{Philippe}\binits{P.}},
\bauthor{\bsnm{Mischler}, \bfnm{St{\'e}phane}\binits{S.}} \AND
\bauthor{\bsnm{Perthame}, \bfnm{Beno{\^{\i}}t}\binits{B.}}
(\byear{2004}).
\btitle{General entropy equations for structured population models and
scattering}.
\bjournal{C. R. Math. Acad. Sci. Paris}
\bvolume{338}
\bpages{697--702}.
\bid{doi={10.1016/j.crma.2004.03.006}, mr={2065377}}
\end{barticle}
%
\endbibitem

\bibitem{bblMichelMischlerPerthame2}
%
\begin{barticle}[mr]
\bauthor{\bsnm{Michel}, \bfnm{Philippe}\binits{P.}},
\bauthor{\bsnm{Mischler}, \bfnm{St{\'e}phane}\binits{S.}} \AND
\bauthor{\bsnm{Perthame}, \bfnm{Beno{\^{\i}}t}\binits{B.}}
(\byear{2005}).
\btitle{General relative entropy inequality: An illustration on growth models}.
\bjournal{J. Math. Pures Appl. (9)}
\bvolume{84}
\bpages{1235--1260}.
\bid{doi={10.1016/j.matpur.2005.04.001}, mr={2162224}}
\end{barticle}
%
\endbibitem

\bibitem{bblOksendal}
%
\begin{bbook}[vtex]
\bauthor{\bsnm{{\O}ksendal}, \bfnm{Bernt}\binits{B.}}
(\byear{2003}).
\btitle{Stochastic Differential Equations},
\bedition{6th} ed.
\bpublisher{Springer}, \baddress{Berlin}.
\bid{mr={2001996}}
\end{bbook}
%
\endbibitem

\bibitem{bblOssiander}
%
\begin{barticle}[mr]
\bauthor{\bsnm{Ossiander}, \bfnm{Mina}\binits{M.}}
(\byear{2005}).
\btitle{A probabilistic representation of solutions of the incompressible
{N}avier--{S}tokes equations in {$\bold R\sp3$}}.
\bjournal{Probab. Theory Related Fields}
\bvolume{133}
\bpages{267--298}.
\bid{doi={10.1007/s00440-004-0418-z}, mr={2198702}}
\end{barticle}
%
\endbibitem

\bibitem{bblRozovskiBook}
%
\begin{bbook}[vtex]
\bauthor{\bsnm{Rozovski{\u\i}}, \bfnm{B. L.}\binits{B. L.}}
(\byear{1990}).
\btitle{Stochastic Evolution Systems: Linear Theory and Applications to Nonlinear Filtering}.
\bseries{Mathematics and Its Applications (Soviet Series)}
\bvolume{35}.
\bpublisher{Kluwer Academic}, \baddress{Dordrecht}.
(\bnote{Translated from
the Russian by A. Yarkho}.)
\bid{mr={1135324}}
\end{bbook}
%
\endbibitem

\bibitem{bblRuban}
%
\begin{barticle}[auto:SpringerTagBib|2010-03-24|17:41:21]
\bauthor{\bsnm{Ruban}, \bfnm{V. P.}\binits{V. P.}}
(\byear{1999}).
\btitle{Motion of magnetic flux lines in magnetohydrodynamics}.
\bjournal{JETP}
\bvolume{89}
\bpages{299--310}.
\end{barticle}
%
\endbibitem

\bibitem{bblThomannOssiander}
%
\begin{bincollection}[mr]
\bauthor{\bsnm{Thomann}, \bfnm{Enrique}\binits{E.}} \AND
\bauthor{\bsnm{Ossiander}, \bfnm{Mina}\binits{M.}}
(\byear{2003}).
\btitle{Stochastic cascades applied to the {N}avier--{S}tokes equations}.
In \bbooktitle{Probabilistic Methods in Fluids}
\bpages{287--297}.
\bpublisher{World Scientific, River Edge, NJ}.
\bid{doi={10.1142/9789812703989_0019}, mr={2083379}}
\end{bincollection}
%
\endbibitem

\bibitem{bblWaymire}
%
\begin{bmisc}[vtex]
\bauthor{\bsnm{Waymire}, \bfnm{Edward C.}\binits{E. C.}} (\byear{2002}).
\bhowpublished{Multiscale and multiplicative processes in fluid flows.
In \textit{Instructional and Research Workshop on Multiplicative Processes and Fluid
Flows (MaPhySto, Aarhus Univ., 2001). Lectures on
Multiscale and Multiplicative Processes in Fluid Flows} \textbf{11}.
Available at
\href{http://www.maphysto.dk/cgi-bin/gp.cgi?publ=407}{http://www.maphysto.dk/cgi-bin/}
\href{http://www.maphysto.dk/cgi-bin/gp.cgi?publ=407}{gp.cgi?publ=407}}.
\end{bmisc}
%
\endbibitem

\bibitem{bblWaymire2}
%
\begin{barticle}[mr]
\bauthor{\bsnm{Waymire}, \bfnm{Edward C.}\binits{E. C.}}
(\byear{2005}).
\btitle{Probability \& incompressible {N}avier--{S}tokes equations: An overview
of some recent developments}.
\bjournal{Probab. Surv.}
\bvolume{2}
\bpages{1--32 (electronic)}.
\bid{doi={10.1214/154957805100000078}, mr={2121794}}
\end{barticle}
%
\endbibitem

\bibitem{bblWebber}
%
\begin{barticle}[auto:SpringerTagBib|2010-03-24|17:41:21]
\bauthor{\bsnm{Webber}, \bfnm{W.}\binits{W.}}
(\byear{1968}).
\btitle{\"Uber eine Transformation der hydrodynamischen Gleichungen}.
\bjournal{J. Reine Angew. Math.}
\bvolume{68}
\bpages{286--292}.
\end{barticle}
%
\endbibitem

\bibitem{bblZhang}
%
\begin{barticle}[auto:SpringerTagBib|2010-03-24|17:41:21]
\bauthor{\bsnm{Zhang}, \bfnm{Xicheng}\binits{X.}}
(\byear{2010}).
\btitle{A stochastic representation for backward incompressible
Navier--Stokes
equations}.
\bjournal{Probab. Theory Related Fields}
\bvolume{148}
\bpages{305--332}.
\bid{mr={2653231}}
%
\end{barticle}
%
\endbibitem

\end{thebibliography}
\end{document}